\documentclass[10pt]{article}
  \usepackage[reqno, namelimits, sumlimits]{amsmath}
        \usepackage{amssymb, amsfonts}
        \usepackage{mathrsfs}
        \usepackage[margin=1.4in, paperwidth=8.5in, paperheight=11.0in]{geometry}

\newcommand{\be}{\begin{equation}}
\newcommand{\ee}{\end{equation}}
\newcommand{\bb}{\bigskip}
\newcommand{\la}{\label}

\newcommand{\rf}[1]{(\ref{#1})}

\newcommand{\ve}{\varepsilon}
\newcommand{\om}{\omega}
\newcommand{\Om}{\Omega}
\newcommand{\vf}{\ensuremath{\varphi}}
\newcommand{\ddt}{\frac{\partial}{\partial\theta}}
\newcommand{\ddz}{\frac{\partial}{\partial z}}

\renewcommand{\Re}{\operatorname{Re}}
\renewcommand{\Im}{\operatorname{Im}}
\newcommand{\ZZ}{\mathbf Z}

\newcommand{\R}{{\mathbf R}}

\newcommand{\Z}{{\mathcal Z}}
\newcommand{\C}{{\mathbf C}}

\newcommand{\CC}{\mathbf C}

\newcommand{\OO}{\mathcal O}
\newcommand{\PP}{\mathcal P}
\newcommand{\HH}{\mathcal H}

\newcommand{\MM}{\mathcal M}

\newcommand{\YY}{\mathcal Y}

\newcommand{\dd}{\Delta}
\renewcommand{\gg}{\mathfrak g}

\newcommand{\so}{{\bf S}^1}
\newcommand{\hoh}{\dot H^{\frac 12}}

\newcommand{\intmii}{\int_{-\infty}^\infty}

\newcommand{\nyz}[1]{||{#1}||_{Y_0}}

\newcommand{\nlt}[1]{||{#1}||_{L^2}}
\newcommand{\tL}{\tilde L}

\newcommand{\htt}{H^{\frac32}}
\newcommand{\tv}{\tilde v}

\newcommand{\sm}{\smallskip}

\newcommand{\pul}{\frac12}
\newcommand{\less}{\lesssim}
\newcommand\myeq{\stackrel{{\normalfont\mbox{\tiny def}}}{=}}

\def\XXint#1#2#3{{\setbox0=\hbox{$#1{#2#3}{\int}$ }
\vcenter{\hbox{$#2#3$ }}\kern-.6\wd0}}

\newcommand{\iso}{\int_{\so}}
\newcommand{\tY}{\tilde Y}
\newcommand{\tYY}{\tilde\YY}
\newcommand{\ymz}{Y^{(m)}_0}
\newcommand{\ym}{Y^{(m)}}

\newcommand{\ed}{\end{document}}

\newcommand{\cara}{\rule{2mm}{2mm}}

\usepackage{amsmath}
\numberwithin{equation}{section}
\newtheorem{theorem}{Theorem}[section]
\newtheorem{lemma}{Lemma}[section]
\newtheorem{corollary}{Corollary}[section]

\begin{document}
\title{On the De Gregorio modification of the Constantin-Lax-Majda model}
\author{H.~Jia%\footnote{IAS Princeton during 2016--2017}
\,,\,\,
 S.~Stewart\,, \,V.~Sverak
 \\ \\
{\small  University of Minnesota}
\\
}

\date{}
\maketitle
\begin{abstract}We study a modification due to De Gregorio of the Constantin-Lax-Majda (CLM) model $\om_t=\om\, H\!\om$ on the unit circle. The De Gregorio equation is \hbox{$\om_t+u\om_x-u_x\om = 0,$} \hbox{$u_x=H\om$}. In contrast with the CLM model, numerical simulations suggest that the solutions of the De~Gregorio model with smooth initial data exist globally for all time, and generically converge to equilibria when $t\to\pm\infty$, in a way resembling inviscid damping. We prove that such a behavior takes place near a manifold of equilibria.
\end{abstract}

\section{Introduction}
The Constantin-Lax-Majda (CLM) model \cite{CLM} is
\be\la{1}
\om_t=\om H\om\,,
\ee
where $H$ is the Hilbert transform. It can be considered on the real line $\R$ or on the circle $\so$. We will mostly work on the circle. We will use the coordinates
\be\la{2}
z=e^{i\theta}\,,
\ee
where $z$ is considered as a complex number and $\theta\in (-\pi,\pi]$.
De Gregorio~\cite{DG} suggested the following modification of~\rf{1}:
\be\la{2b}
\om_t+u\om_\theta=u_\theta \om\,,\qquad u_\theta = H\om\,,\quad \int_{\so} u=0\,.
\ee
If we consider $\om, u$ as vector fields on $\so$, we can write~\rf{2b} as
\be\la{3}
\om_t+[u,\om]=0\,,
\ee
where $[a,b]=a_\theta b-a b_\theta$ is the usual Lie bracket for vector fields.

As observed by Constantin, Lax, and Majda, the CML model is equivalent to the local equation for holomorphic functions on the unit disc $D$
\be\la{x1}
\dot f=-if^2\,,
\ee
where $f$ is the holomorhic extension of the function $\om+iH\om$ into $D$. From this it is clear that the solution with initial datum $\om_0$ blows up in finite (positive) time if and only if the image of $\so$ by the function $f_0=\om_0+i H\om_0$ intersects the positive imaginary axis $(\R_+)\,i$.
From~\rf{x1} it is clear that the CLM model~\rf{1} is locally-in-time well-posed for initial conditions $\om_0$ satisfying $||\om_0||_{L^\infty}+||H\om_0||_{L^\infty}<\infty$. In particular, it is locally-in-time well-posed for $\om_0\in H^s(\so)$ for $s>\frac12$. On the other hand, it is also clear that it is not locally-in-time well-posed for a general $\om_0\in C(\so)$ (continuous functions).\footnote{ By contrast, it is natural to expect that the transport equation $\om_t+u\om_x=0$, with the Biot-Savart law as in~\rf{2b} is locally-in-time well-posed in $C(\so)$.}
Based on the comparison with the CLM-model, is natural to expect that the De Gregorio equation~\rf{2b} is locally-in-time well-posed in $H^s(\so)$ for $s>\frac 12$, but not in $H^{\frac 12}(\so)$ or $C(\so)$. One also expects the Beale-Kato-Majda-type criterion: the $L^\infty_tH^s_x$-regularity of the solution (for $s>\frac 12$)  in a closed interval $[0,T]$ should be controlled by the condition
\be\la{BKMi}
\int_0^T||\om(t)||_{L^\infty}\,dt<+\infty\,,
\ee and, in particular, when the integral in~\rf{BKMi} is finite, the solution can be locally continued beyond $T$ without a loss of regularity.  This is proved for $s\ge 1$ in~\cite{BKP}.  More discussion of these topics is in subsection~\ref{prelim6}.
Very recently, finite-time blow-up for initial data in low regularity spaces (allowing infinite derivatives) in which the equation is still locally well-posed was proved in~\cite{EJ}.

Numerical simulations seem to suggest that there is no blow-up from smooth initial data for the De Gregorio equation, as already observed in~\cite{OSW}.  Our own numerical simulations suggest that for generic smooth initial data the solutions $\om(t)$ approach equilibria $A\sin(\theta-\theta_0)$ for $t\to \pm\infty$, although not in smooth norms. The convergence is only in $H^s$ for $s<\frac 32$, and not in $H^{\frac 32}$. Moreover, the initial datum has to be sufficiently regular, with $\om_0\in H^{\frac 32 +\epsilon}$ for any $\ve>0$ being probably sufficient, but $\om_0\in H^{\frac 32}$ presumably allowing much more complex dynamical behavior.
A good toy model for these phenomena is the linear equation
\be\la{x2}
\om_t+[b,\om]=0\,,
\ee
with $b=\sin\theta$, which can be completely analyzed by explicit calculation.

Regularising effects of transport terms have been observed in other models, see for example~\cite{HouLi}.

The only known conserved quantities for the De Gregorio equation are the orbit invariants discussed in subsection~\ref{orbits}, and the quantity $\iso \om(\theta)\,d\theta$. None of these are coercive. The conjectured long-time behavior, together with the orbit invariants, seems to put strong constraints on possible conserved quantities, and it is not clear if there is any good coercive conserved quantity at all. This should be contrasted with the remarkably good numerical behavior of the solutions, and their apparent convergence to steady states for $t\to \pm\infty$ for smooth data.

We will study the dynamics near the equilibria $A\sin(\theta-\theta_0)$. Our main theoretical result is the following.
\begin{theorem}\label{mainresult} {\rm (Non-linear stability of equilibria)}
Let $\om_0$ be a $C^2$ function which is $C^2$-close to an equilibrium $\Om_{A,\theta_0}=A\sin(\theta-\theta_0)$. Then the De Gregorio equation~\rf{2b} has a unique $C^2$-solution $\om(t)$ with $\om(0)=\om_0$ defined for all $t\in\R$. Moreover, as $t\to\pm\infty$, the solution $\om(t)$ approaches equilibria $\Om_{A^{\pm},\theta_0^{\pm}}$, respectively, for suitable $A^\pm$ and $\theta_0^\pm$. The convergence is exponential in $H^{s}$ for any $s<\frac 32$, but there is no strong convergence in $H^\frac 32$.
The amplitudes $A^{\pm}$ can be determined explicitly from $\om_0$ as described in Remark 1 below.
\end{theorem}

As we will see, the higher norms $||\om(t)-\Om_{A^\pm,\theta_0^\pm}||_{H^s}$ with $s>\frac 32$ typically grow exponentially as $t\to\pm\infty$.

\medskip\noindent
{\bf Remark 1.} The proof shows that the amplitudes $A^\pm$ are determined as follows. We first note that when $\om_0$ is sufficiently close to $\Om_{A,\theta_0}$ in $C^2$, then $\om_0$ has exactly two zeroes. At one of them the derivative $\om_{0x}$ is positive and it is negative at the other one.  Let us denote the former one by $x_1$ and the latter one by $x_2$. Then $A^+=-\om_{0x}(x_2)$ and $A^-=\om_{0x}(x_1)$.

\medskip

The proof of the theorem is based on a careful analysis of the linearized operator, in suitable moving frames. Crucial points of the proof include establishing an (almost) unitarity of the linearized evolution in $H^{\frac32}$, ruling out the point spectrum in $H^{\frac 32}$ by the use of ODEs in the complex domain, establishing connection of the linearized equation with the simple evolution~\rf{x2}, establishing exponential decay of linearized solutions in an auxiliary space $Y_0$, and finally using bootstrapping to handle the non-linearity.

Other aspects of the De Gregorio equation and its modifications are studied for example in~\cite{CastroCordoba, BKP, EK, EKW, Wunsch}. Some of these references discuss geometric aspects of the equation, in the spirit of Arnold-Khesin~\cite{ArnoldKhesin}, although one has to replace Levi-Civita connections of a Riemannian metric with more general connections, as discussed in~\cite{EKW}.

\section{Preliminaries}\label{prelim}
\subsection{Simple observations}
Denoting by $\phi^t$ the diffeomorphism of the circle defined by the flow
\be\la{4}
\dot x(t)=u(x(t),t)
\ee
with $\phi^0={\rm id}$, equation~\rf{3} is equivalent to
\be\la{5}
\om(t)=\phi^t_{\#}\om(0)\,,
\ee
where $\phi_{\#}a$ is the push-forward of a vector field $a$ by the diffeomorphism $\phi$, i.\ e.\
\be\la{6}
\phi_{\#}a\,(\theta)=\phi'(\phi^{-1}(\theta))\,a(\phi^{-1}(\theta))\,.
\ee
When $u$ and $\om$ are independent, equation~\rf{3} is invariant under diffeomorphisms:
\be\la{8}
(\phi_{\#}\om)_t+[\phi_{\#}u,\phi_{\#}\om]=0
\ee
for any diffeomorphism $\phi$ of $\so$. However, the ``Biot-Savart law"
\be\la{9}
u_\theta=H\om
\ee
is only invariant under a much smaller group of transformations, and hence the non-linear equation does not seem to have any (exact) symmetries beyond the obvious ones given by rotations and reflections.

It is perhaps worth noting that the operator
\be\la{10}
\Lambda = -H \partial_\theta
\ee
considered on {\it scalar functions} (as opposed to vector fields) on the circle $\so$ is covariant under the projective transformations of the circle in a similar way as $-\dd$ on the scalar functions on the disc is invariant under the conformal transformations of the disc.
This can be easily seen for example as follows. The quadratic form
\be\la{11}
(\om,\om)_{\dot H^\pul}=\iso (-\Lambda \om) \om\,dx
\ee
can also be expressed as
\be\la{12}
(\om,\om)_{\dot H^\pul}=\int_{D} |\nabla h|^2\,dx\,dy\,,
\ee
where $h$ is a harmonic extension of $\om$ from the circle $\so$ to the unit disc $D$.
Now the conformal diffeomorphisms $\gamma\colon D\to D$ leave the Dirichlet integral on the right-hand side of the last equation invariant, if we act on $h$ by $h\to h\circ \gamma^{-1}$. The restrictions of all possible $\gamma$ to $\so$ give exactly the orientation-preserving pojective transformations of $\so$. This implies
\be\la{13}
\Lambda (\om\circ\gamma)=|\gamma'| \!\cdot\,(\Lambda\om)\circ\gamma\,.
\ee

\subsection{Orbit invariants}\label{orbits}
Let $G$ be a Lie group and let $\gg$ be its Lie algebra. Consider the equation for a $\gg-$valued function of time
\be\la{o1}
\dot \xi=[L(\xi),\xi]\,,
\ee
where $L\colon\gg\to\gg$ is a smooth function. (The equation is in the Lax form.) For our purpose here we can think of the case when $L$ is linear. The trajectory $\xi(t)$ with $\xi(t_0)=\xi_0$ lies in the adjoint orbit
\be\la{o2}
\OO_{\xi_0}=\{a\cdot\xi_0\cdot a^{-1}\,,\,a\in G\}\,,
\ee
and hence the orbit invariants, such as eigenvalues of (a suitable representation of) $\xi$ are preserved.
The De Gregorio equation is formally of this form, with $G={\rm Diff}_+(\so)$, the group of the orientation-preserving diffeomorphisms of the circle, except that ${\rm Diff}_+\,(\so)$ is infinite-dimensional,  and hence some of the properties of finite-dimensional Lie groups may not be available.
The adjoint orbits in this case are
\be\la{i3}
\OO_{\om_0} =\big\{\phi_{\#}\om_0\,,\,\phi\in{\rm Diff}_+\,(\so)\big\}\,.
\ee
What are the invariants of such orbits? In general this is not an easy question. In the case when $\om_0$ has finitely many zeroes of finite order, the full classification was obtained by Hitchin~\cite{Hitchin}. Here we will only discuss the case when the zeroes of $\om_0$ are non-degenerate in the sense that $\om'_0(\theta)\ne 0$ when $\om_0(\theta)=0$.
Let $\theta_1<\theta_2<\dots<\theta_{2m}$ be such zeroes, and let
\be\la{i4}
(a_1,\dots,a_{2m})=(\om'_0(\theta_1),\dots,\om'_{0}(\theta_{2m}))\,.
\ee
Also, let us define
\be\la{i5}
b={\rm p. \, v.}\,\iso \frac {d\theta}{\om_0(\theta)}\,,
\ee
where p. v. means the principle value.
Then the data
\be\la{i6}
2m, (a_1,\dots, a_{2m}), b\,
\ee
where $(a_1,\dots, a_{2m})$ is considered modulo  cyclic permutations of $a_1,\dots,a_{2m}$, are invariants of the orbit. In our non-degenerate case this list of invariants is complete, i.\ e.\ two orbits with the same invariants coincide, see~\cite{Hitchin}.

The conservation of the derivative $\om_\theta(\theta(t),t)$ at the zeroes of $\om$ is seen easily directly from the De Gregorio equation~\rf{2b}.
Taking a derivative of the equation, we obtain
\be\la{i7}
\om_{\theta t}+u\omega_{\theta\theta}=u_{\theta\theta}\om\,.
\ee
If $\om(\theta(t),t)=0$, then~\rf{i7} implies that $\om_\theta(\theta(t),t)$ is preserved as $t$ changes. Clearly this remains true for any $u(x,t)$ (it does not have to be given by a specific Biot-Savart law), reflecting again the fact that $\om_\theta(\theta(t),t)$ at $\theta(t)$ with $\om(\theta(t),t)=0$ is an invariant of the orbit.
This ``conservation law", together with the conservation of ${\rm p.\,v.\,} \iso d\theta/\om(\theta,t)\,,$ is an analogue of the Kelvin-Helmholtz law for the classical fluids.

\subsection{Conservation of $\iso \om\,d\theta$}

In general the integral $\iso\om(\theta)\,d\theta$ is not invariant on the adjoint orbit, but for the evolution by equation~\rf{2b} it is invariant, as we have
\be\la{i8}
\frac{d}{d t}\iso\om(\theta,t)\,d\theta=\iso (-u\om_\theta+u_\theta\om)=\iso 2u_\theta\om=\iso 2(H\om) \om = 0\,.
\ee
In what follows we will work with the solutions $\om$  satisfying
\be\la{i9}
\iso \om(\theta,t)\,d\theta = 0\,.
\ee
In the general case
\be\la{i9a}
\iso \tilde\om(\theta,t)\,d\theta=c
\ee
we can can replace $\tilde\om$ by $\om+c$, where $\om$ still satisfies~\rf{i9}. The equation then becomes
\be\la{i9b}
\om_t+[u,\om]=cH\om\,.
\ee
The solutions of~\rf{i9b} corresponding to the steady states $A\sin(\theta-\theta_0)$ of \rf{3} become
\be\la{i9c}
\om(x,t)=A\sin(\theta-\theta_0-ct)\,.
\ee

\subsection{Other choices of gauge for the velocity field}
Let $\theta_0\in(-\pi,\pi]$.
It will be sometimes useful for us to change coordinates and instead of ``calibrating" the velocity field $u$ by
\be\la{i9e}
\iso u \,d\theta = 0\,,
\ee
we will modify by it by a constant (depending on time) and work with the field $\tilde u(x,t)$ defined by
\be\la{i9f}
\tilde u_\theta = H\om\,,\qquad \tilde u(\theta_0,t)=0\,.
\ee
Assume that $\om_1(\theta,t)$ is a solution of~\rf{2b}, with the corresponding vector field $u_1$, and set
\be\la{i9g}
\om(\theta,t)=\om_1(\theta+\vartheta(t),t)\,,\qquad u(\theta,t)=u_1(\theta+\vartheta(t),t)\,,
\ee
where $\vartheta$ is a function of time.
Then
\be\la{i9h}
\om_{t}+[u-\dot\vartheta(t)\,, \,\om]=0.
\ee
We see that if we choose $\vartheta$ so that
\be\la{i9i}
\dot\vartheta(t)=u(\theta_0,t)\,,
\ee
which amounts to solving $\dot \vartheta(t)=u_1(\theta_0+\vartheta(t),t)$,
the field $\om$ will solve
\be\la{i9j}
\om_t+[\tilde u,\om]=0\,.
\ee
If we start with a solution of~\rf{i9j}, we can obtain a solution of~\rf{3} by a similar change of variable. We see that the equations~\rf{i9j} and~\rf{3} are equivalent. When more convenient, we will work with~\rf{i9j} rather than~\rf{3}.

\subsection{Equilibria and numerically observed long-time behavior}
It is easy to see that functions of the form
\be\la{e1}
\om(\theta)= A\sin m(\theta-\theta_0)\,,\quad A\in\R\,,\,\,\theta_0\in(-\pi,\pi]\,,\quad m\in\ZZ
\ee
are equilibria of~\rf{3}. For each $m$ these form a two-dimensional manifold $\MM_m$ in the space of smooth (real-valued) functions on the circle $\so$.
Numerical experiments (performed in Matlab using a simple pseudo-spectral method) suggest that the manifolds $\MM_m$ with $m\ge 2$ are unstable, and that generic smooth solutions of~\rf{3} approach the manifold $\MM_1$, although not in the space of smooth (or even $C^1$) functions.  The convergence to equilibria appears to hold in $\dot H^s$ for $s<3/2$. This is consistent with the conservation of the orbits $\OO_{\om_0}$ and the invariants~\rf{i6}. We note that, in general,  a $C^1$ convergence (or $\dot H^{\frac 3 2}$ convergence) would not be consistent with the conservation of these quantities. This is obvious for the $C^1$ case;
the case of $\dot H^{\frac 32}$ follows from the analysis below.
In fact, a generic trajectory $\om(t)$ appears to have  well-defined limits $\om_\infty$ and $\om_{-\infty}$ as $t\to +\infty$ and $t\to-\infty$ respectively (in topologies just below $\dot H^{\frac 32}$ regularity).

Previous numerical results were reported for example in~\cite{OSW}. The results there agree with our numerical observation that there appears to be no blow-up.
\sm
The approach to equilibria seems to happen in a way which is similar to the following linear toy model. Let
\be\la{M2}
b(\theta)\frac{\partial}{\partial \theta}
\ee
be a smooth vector field on $\so$ with exactly two non-degenerate equilibria at $\theta=\theta_1$ and $\theta=\theta_2$ on $S^1$.
We can assume the equilibrium at $\theta_1$ is unstable, i.\ e.\ $b'(\theta_1)>0$ (together with $b(\theta_1)=0$). Then necessarily $b'(\theta_2)<0$. (Recall that we assume the equilibria are non-degerate.)
Consider now the equation
\be\la{M3}
\om_t+[b,\om]=0
\ee
with the initial value $\om(\theta,t)|_{t=0}=\om_0(\theta)$, where $\om_0$ is a smooth initial condition, where $\om_0(\theta_1)=0\,,\,\om_0'(\theta_1)=\alpha\,.$
In this situation it is not hard to show that
\be\la{M4}
\lim_{t\to\infty} \om(\theta,t)=\alpha b(\theta)\,,\quad \theta\ne \theta_2\,,
\ee
and the convergence is uniform on compact subsets of $\so\setminus\{\theta_2\}$. One way to see this is to change the coordinate $\theta$ to $\xi$ so that the vector field $b$ in $\so\setminus\{\theta_2\}$ becomes
\be\la{M5}
c\xi\frac{\partial}{\partial\xi}\,\quad \xi\in \R\,.
\ee
The flow map $\phi^t$ given by the last field is given explicitly by
\be\la{M6}
\phi^t(\xi)=e^{ct}\xi
\ee
and the limit $\phi^t_{\#}\om_0$ for $t\to\infty$ is easily calculated, after expressing $\om_0$ in the coordinate $\xi$ via
\be\la{M7}
\om_0(\theta)\frac{\partial}{\partial\theta}=\frac {\om_0(\theta)}{b(\theta)}\, b(\theta)\frac{\partial}{\partial\theta}=
\frac{\om_0}{b}\,c\xi\frac{\partial}{\partial\xi}\,.
\ee

In the linear example~\rf{M3} we had to assume that the zeroes of $b$ and $\om_0$ were ``aligned". (If $\om_0(\theta_1)\ne0$, it is easy to see that $|\om(\theta,t)|\to\infty$ for $\theta\ne\theta_2$.)
The non-linear equation seems to be able to align the zeroes of $\om$ and $u$ ``by itself".

\subsection{Local well-posedness for  $\om_0\in H^{\pul+\epsilon}$ and the  BKM criterion}\label{prelim6}
The local-in-time well-posedness for $\om_0\in H^1$ is proved in~\cite{OSW}, and the Beale-Kato-Majda-type criterion mentioned in the introduction, namely that the $L^\infty_t H^s_x$ regularity in any closed interval $[0,T]$ is controlled by the condition
$$
\int_0^T||\om(t)||_{L^\infty}\,dt < +\infty\,,
$$
 is proved in~\cite{BKP} when $s\ge 1$. One can generalize these results to $s>\frac 12$ based on the methods of~\cite{Danchin, Chemin-book}, and also~\cite{Ebin-Marsden, BKP, Lee-Preston}. Here we only briefly outline the arguments, leaving a more detailed exposition of these topics for a future work.

Motivated by the works~\cite{Danchin, Chemin-book, Ebin-Marsden, BKP, Lee-Preston}, we  can re-write the De Gregorio equation~\rf{2b} in terms of $u$ as follows:
\be\la{dg-u-form}
u_t+uu_x=uu_x-\Lambda^{-1}\left( u\Lambda u_x - u_x\Lambda u\right)\myeq B(u,u)\,.
\ee
The operator $\Lambda$ has a one-dimensional kernel consisting of constant functions, but if we work with functions of zero average, we do not have difficulties with the invertibility of $\Lambda$, if we take into account~\rf{i8}. We fix $s\in(\frac 12, 1)$.
The main point now is that the expression $B(u,u)$ on the right-hand side of~\rf{dg-u-form} is a continuous quadratic mapping from $H^{1+s}$ to itself. We notice that, due to cancellations in the expression for $B$, one can expect
\be\la{B-est}
||B(u,u)||_{H^{1+s}}\lesssim ||u_x||_{L^\infty}||u||_{H^{1+s}}\,,
\ee
from which it is not hard to get the Beale-Kato-Majda-type criterion at the level of $L^\infty_tH^s_x$ regularity for $\om$, with the help of the Kato-Ponce commutator estimate~\cite{KP} and an often-used trick of Kato involving the estimation of $u_x$ through $\om$ and $\log$ of a higher norm, see for example Proposition 2.104 in~\cite{Chemin-book}.

For the local well-posedness for $\om_0\in H^s$ one can follow (with some modifications) either the methods of~\cite{Chemin-book, Danchin}, working with the Eulerian formulation\footnote{and showing the ``quasi-linear wellposedness" in the sense of~\cite{Tzvetkov}},  or~\cite{Ebin-Marsden, Lee-Preston}, working with the Lagrangian formulation, and showing (again with some modifications) that the vector field which defines the equations in Lagrangian coordinates on the tangent space of the group of $H^{1+s}$-diffeomorphisms of $\so$ is Lipschitz, so we are dealing with an ODE in these coordinates, and standard ODE theorems can be applied.

\section{Linearized Stability}
In this section we will study the linearization of equation~\rf{2b} about the steady solution
\be\la{l1}
\Om(\theta)=-\sin\theta\,,
\ee
The corresponding velocity field given by the Biot-Savart law in~\rf{2b} is
\be\la{l2}
U(\theta)=\sin\theta\,.
\ee
The linearized equation is
\be\la{l3}
\eta_t+[U,\eta] + [v,\Om]=0\,,\quad v_\theta=H\eta\,,\quad \iso v=0\,.
\ee
Note that if the equation did not have  term $[v,\Om]$, we would be in the situation of~\rf{x2}, and the long-time behavior of $\eta$ would be easy to determine.

The linearized operator
\be\la{l4}
L(\eta)=-[U,\eta]-[v,\Om]=-[U,\eta+v]\,.
\ee
has two important properties which will help us to handle the situation.

\begin{lemma}\label{commutation}
$L$ commutes with the Hilbert transform $H$.
\end{lemma}

\noindent
{\bf Proof}\\
Recall that we assume that $\iso \eta=0$, and from the Biot Savart law we then see that the three Fourier coefficients of the function $\eta+v$ corresponding to $e^{-i\theta},1$ and $e^{i\theta}$ vanish. This easily gives the result, as commutation by $U$ shift the Fourier frequencies at most by $1$, and the Fourier multiplier of $H$ is constant on the positive frequencies and is also constant on the negative frequencies.\,\,\,\cara

\bb
We note that Lemma~\rf{commutation} requires that the functions $\eta$ satisfy $\eta_0=0$. If we wish to work with the natural extension of $L$ to $\eta_0\ne 0$ defined by $L e_0=\pul e_{-1}+\pul e_1$, the commutator $[L,H]$ will not vanish on $e_0$. However, it will still vanish if we mod out by the linear span of $e_{-1},e_0,e_1$ which is a subspace invariant  under both $H$ and (the extended) $L$.
\bb

We now aim to show that $L$ is skew-symmetric with respect to a certain quadratic form. This calculation seems to be easiest in the Fourier variables. For $k\in\ZZ$ we denote by $e_k$ the function $e^{ik\theta}$.
A direct calculation shows that
\be\la{l5}
Le_k=A_ke_{k-1}+B_k e_{k+1}\,,\quad k\ne 0\,,
\ee
where, for $k\ne 0$,
\be\la{l6}
A_k=\pul(k+1)(1-\frac 1{|k|})\,,\qquad B_k=\pul(-k+1)(1-\frac1{|k|})\,.
\ee
In terms of the Fourier coordinates $\eta_k$ this means that
\be\la{l7}
(L\eta)_k=B_{k-1}\eta_{k-1}+A_{k+1}\eta_{k+1}\,,\qquad k\ne 0\,,
\ee
where we adopt the convention $B_0\eta_0=0$ and $A_0\eta_0=0$. (We recall that we work with functions $\eta$ for which $\eta_0=0$. A natural extension of the operator to functions with $\eta_0\ne 0$ is by setting $A_0=B_0=\pul$, as the term $1-\frac 1{|k|}$ in~\rf{l5} arises as $1-\frac{\rm sign k}{k}$.)
We note that $A_k=0\,,\,\,B_k=0$ for $|k|=1$, and that $A_{-k}=B_k, B_{-k}=A_k$.

The formulae~\rf{l6} also provides a proof of Lemma~\ref{commutation}.

The evolution equations for $\eta_k\,,\,\,k=1,2,3,\dots$  are
\be\la{l8}
\begin{array}{rcl}
\dot \eta_1 & = & A_2\eta_2\\
\dot \eta_2 & = & A_3\eta_3\\
\dot \eta_3 & = & B_2\eta_2+A_4\eta_4\\
 &\dots &
\end{array}
\ee
Note that the system for $\eta_2,\eta_3,\dots$ is closed and the variable $\eta_1$ can be calculated by integration of $A_2\eta_2$ once the components $\eta_2,\eta_3,\dots$ are known.

We aim to find  $c_k>0\,,\,k=2,3,\dots$ so that $L$ is anti-hermitian with respect the hermitian form
\be\la{l9}
\langle \eta,\eta\rangle = \sum_{k=2}^\infty c_k\eta_k\bar\eta_k\,.
\ee
An easy calculation shows that the condition on the sequence $c_k$ is
\be\la{l10}
c_{k+1}=\frac{A_{k+1}}{-B_k}\,c_k\,,\quad k=2,3,\dots
\ee
and hence we can write
\be\la{l11}
c_{k+1}=A_{k+1}A_{k}A_{k-1}\frac {A_{k-2}A_{k-3}\dots A_3}{(-B_k)(-B_{k-1})\dots (-B_2)}\,\,c_2\,.
\ee
It is easy to see that the fraction on the right-hand side of the last equation has a finite strictly positive limit as $k\to\infty$ and hence we see that
\be\la{l12}
c_k\sim k^3\,.
\ee
In fact, a more detailed calculation shows
\be\la{ckprecise}
c_k=(k-1)^2(k+1)\,.
\ee
Also, the conservation of the form~\rf{l9} by the evolution  given by $\dot\eta=L\eta$ can be seen directly by formulating the evolution in a Hamiltonian form. On the phase-space given by the (complex) coordinates
$\eta_1,\eta_2,\dots$ we define the Hamiltonian
\be\la{Hamiltonian}
{\mathscr H}={\mathscr H}(\eta)=\sum_{k=1}^\infty \frac14 c_k \eta_k\bar\eta_k
\ee
Consider the (infinite) matrix
\be\la{matrixJ}
J=\left(
\begin{array}{rrrrr}
0 & a_1 & 0 & 0 & \dots\\
-a_1 & 0 & a_2 & 0 & \dots\\
0 & -a_2 & 0 & a_3 & \dots\\
0 & 0 & -a_3& 0 &\dots\\
\dots &\dots & \dots& \dots& \dots
\end{array}
\right)\,, \qquad a_k=\frac 1{k(k+1)}\,.
\ee
We can then write our linear system~\rf{l4} in the Hamiltonian from
\be\la{Hfrom}
\dot \eta = L\eta= J \,D{\mathscr H} (\eta)\,,
\ee
which transparently shows that ${\mathscr H}$ is preserved by the evolution.

Note that on the Fourier side $L$ commutes with complex conjugation:
\be\la{even-odd-sym}
L{\bar \eta}=\overline {L\eta}\,, \qquad \eta=(\eta_1,\eta_2,\dots)\,.
\ee
This reflects the fact that in the physical space the evolution preserves the spaces of odd and even functions, respectively.

We see that a good space in which we can consider our equation is the space
\be\la{l13}
X=\dot\HH^{\frac 32}/\CC e_1\,,
\ee
where $\dot\HH^{\frac 32}$ is the space of analytic function on the unit disc with the restriction to the boundary belonging to the Sobolev space $\dot H^{\frac 32}$, and $\CC e_1$ is the subspace of functions which are a multiple of~$z$. Our convention here and in  other similar situations is that $\dot\HH^{\frac32}$ is (equivalent to) a factor space $\HH^{\frac32}/\CC e_0$, so that $X$ is equivalent to $\HH^{\frac32}/(\CC e_0+\CC e_1)$. The variables $\eta_2,\eta_3,\dots$ can be used as coordinates in this space and the hermitian product will be taken as~\rf{l9}.

The  fact that we can restrict our attention to analytic functions can be seen directly from Lemma~\ref{commutation}: we use standard decomposition
\be\la{l14}
\eta=\eta^++\eta^-\,,\quad \eta^+=\pul(\eta+iH\eta)\,,\quad\eta^-=\pul(\eta-iH\eta)\,,
\ee
and due to Lemma~\ref{commutation} we can deal with $\eta^+$ and $\eta^-$ separately.

In what follows we will work with the holomorphic part of $\eta$ and will slightly abuse notation by assuming that $\eta$ is a holomorphic function, and the same for $v$. For holomorphic functions  on the disc it is natural to use the variable $z=e^{i\theta}$. In the holomorphic situation the Biot-Savart law is given by a differential operator:
\be\la{l15}
\eta(z)=-zv'(z), \qquad v'=\frac {dv}{dz}\,.
\ee
To write the equations in the $z-$variable, we use
\be\la{l16}
\frac{\partial}{\partial \theta}=iz\frac{\partial}{\partial z}\,
\ee
and write
\be\la{l17}
\sin\theta \ddt=\pul(z^2-1)\ddz\,.
\ee
The evolution equation for holomorphic $\eta$ can be written symbolically as
\be\la{l18}
\eta_t iz\ddz+\left[\pul(z^2-1)\ddz\,,\, (\eta+v)iz\ddz\right]=0\,,
\ee
which is the same as
\be\la{l19}
\eta_{t} z+\pul(z^2-1)(\eta z + v z)'-z^2(\eta+v)=0\,,\quad \eta=-zv'\,.
\ee
The evolution given by this equation will be unitary in $X$ (as the operator $L$ is anti-hermitian). In addition, the vector $e_2\in X$ (where we slightly abuse notation by using $e_2$ for what really is the projection of $e_2$ into $X$) is a {\it totalizer} for the operator $L$, in the sense that the vectors $e_2, Le_2, L^2e_2, \dots$  generate a dense subspace of $X$. Spectral theory for unitary semigroups now implies that the spectrum of $L$ is purely imaginary, and that there is a measure $\mu$ and an isometry $T\colon X\to L^2(\R,\mu)$ such that
$\tilde L = TLT^{-1}$ is given by
\be\la{l20}
(\tilde L f)(s)=is f(s)\,,\qquad f \in L^2(\R,\mu)\,.
\ee
The corresponding evolution equation in $L^2(\R,\mu)$, equivalent to the projection of~\rf{l18} to $X$, then is
\be\la{l20a}
f_t=\tilde L f
\ee
and its solutions are given by
\be\la{l20b}
f(s,t)=f(s,0)e^{ist}\,.
\ee

Our goal is to obtain information about $\mu$. In particular, we would like to show that $\mu$ is absolutely continuous with respect to the Lebesgue measure. That is enough to conclude from~\rf{l20b} (essentially via the Riemann-Lebesgue lemma) that
\be\la{l21}
\eta(t)\rightharpoonup 0 \,\quad \hbox{in $X$ (weak convergence) as $t\to\infty$}\,,
\ee
for any $\eta(0)\in H^{\frac 32}$. Here we slightly abuse the notation by using $\eta(t)$ also for the projection of $\eta(t)$ to $X$. (Therefore the statement says nothing about the first Fourier mode.)
This result is essentially sharp: we will see that for any given $T>0$ there are non-trivial periodic solutions of period $T$ with $\eta(0)$ just missing $H^{\frac 32}$. Also, even with $\eta(0)\in H^{\frac 32}$ the first Fourier mode may in general not have a limit as $t\to\infty$, while at the same time the projection of $\eta(t)$ approach $0$ weakly in $X$.

\subsection{A model problem}\label{model}
In this subsection we look at a simplified model of~\rf{l19}, which already captures its main features and can be solved explicitly. Some of the calculations will also be important for the analysis of~\rf{l19}.

Let us consider the equation
\be\la{m1}
f_t+\sin\theta\, f_\theta=0\,
\ee
on the unit circle. Denote
\be\la{m2}
Mf=-\sin\theta \,f_\theta\,.
\ee
The analogue of~\rf{l5} for $M$ is
\be\la{m3}
Me_k=A_ke_{k-1}+B_ke_{k+1}\,,\quad k\in \ZZ\,,
\ee
where
\be\la{m4}
A_k=\pul k\,,\qquad B_k=-\pul k\,.
\ee
In terms of the Fourier coordinates $f_k$ the equation gives a closed system for $f_1,f_2,\dots$,  and a closed system for $f_{-1},f_{-2},\dots$. The system for $f_1,f_2,\dots$ is
\be\la{m5}
\begin{array}{rcl}
\dot f_{1} & = & A_2f_2\\
\dot f_{2} & = & B_1 f_1+ A_2 f_3\\
\dot f_3   & = & B_2 f_2+A_4 f_4\\
 &\dots &
\end{array}
\ee
and an analogous system holds for $f_{-1}, f_{-2},\dots$. The equation for $f_0$ is
\be\la{m6}
\dot f_0=B_{-1}f_{-1}+A_1 f_1=\pul \left(f_{-1}+f_1\right)\,.
\ee
One checks easily that $M$ is anti-hermitian with respect to the  $\dot H^{\pul}$ (semi-)norm, given by
\be\la{m7}
||f||_{\hoh}=\frac1{2\pi}\sum_k |k||f_k|^2\,.
\ee
and hence the evolution operator given by~\rf{m1} is unitary in $\hoh$.
The Hilbert transform $H$ does not commute with $M$, but does so modulo constants. More precisely, we have
\be\la{m8}
[H,M]f=\frac{i} 2 \left(f_1-f_{-1}\right)e_0\,,
\ee
so that, recalling our convention that constants are factored out in $\dot H^{\pul}$, we can write
\be\la{m9}
[H,M]=0 \quad \hbox{in $\hoh$.}
\ee
The fact that the evolution by $M$ is an isometry on $\hoh$ is also easily seen from~\rf{l17}: the field $\sin\theta \ddt$ on $\so$ has a holomorphic extension to the unit disc, given by~\rf{l17}, and the evolution given by the extension on harmonic function is an isometry of $\dot \HH^{\pul}$, because conformal transformations preserve the Dirichlet integral on the disc.

We  map the unit disc $D=\{z\,,|z|<1\}$ onto the strip $\OO=\{w\,,-\pi/2<\Im w< \pi/2\}$
 via
\be\la{m10}
z\to w=\log\frac{1-z}{1+z}\,,
\ee
where we take the branch of the log function defined by $\log (re^{i\theta})=\log r +i\theta$ for $r\in (0,\infty)$ and $\theta\in(-\pi,\pi)$. It is easy to check that under the mapping $z\to w$ the vector field
$\pul(z^2-1)\ddz$ is mapped to $\frac{\partial\,}{\partial w}$.
This means that the evolution of the holomorphic sector of $\hoh$ given by~\rf{m1} is equivalent to the equation
\be\la{m12}
f_t+ f_w=0
\ee
in the space $\dot \HH^{\pul}(\OO)$ of holomorphic functions in $\OO$, where the norm on $\dot \HH^\pul(\OO)$ is given by $\int_\OO|f'(w)|^2 \frac{i}{2}\,dw\wedge d\bar w\,$.
The evolution given by ~\rf{m12} is of course $f(w)\to f(w-t)$ and it is diagonalized in the Fourier representation
\be\la{m13n}
f(w)=\int_{-\infty}^{\infty} \vf(s)e^{isw}\,ds\,.
\ee
Assuming $\vf$ is smooth and compactly supported, and using Parserval's identity for the integral over $w_1$, we have
\be\la{m14n}
f'(w)=\int_{-\infty}^{\infty}is\vf(s)e^{isw}\,ds\,.
\ee
 Writing $w=w_1+iw_2$ with $w_1\in(-\infty,\infty)$ and $w_2\in(-\frac\pi 2,\frac \pi 2)$, we have
\be\la{m15n}
\begin{split}
\int_\OO |f'(w)|^2\,dw_1\,dw_2 & =\int_{-\frac\pi 2}^{\frac \pi 2}\,\left(\int_{-\infty}^{\infty}  |f'(w_1+iw_2)|^2\,dw_1\right)\,d w_2\\
& = \int_{-\frac\pi 2}^{\frac \pi 2}\left( 2\pi \int_{-\infty}^{\infty}s^2|\vf(s)|^2e^{-2sw_2}\,ds  \right)\,dw_2\\
& = 2\pi \,\int_{-\infty}^{\infty}|\vf(s)|^2s \sinh\pi s\,\,ds\,.
\end{split}
\ee
Going back to the $z-$variable in the unit disc and remembering $e^w=\frac{1-z}{1+z}$, we see that in the holomorphic sector the spectral decomposition induced by the operator $M$, or equivalently, $\pul(z^2-1)\frac{\partial}{\partial z}$, is given by
\be\la{m16n}
f(z)=\int_{-\infty}^{\infty}\vf(s)\left(\frac{1-z}{1+z}\right)^{is}\,ds\,,
\ee
with
\be\la{m17n}
||f||^2_{\dot \HH{\pul}(D)}\sim \int_{-\infty}^{\infty} |\vf(s)|^2s\,\sinh\pi s\,\,ds\,.
\ee

A simple corollary of the above considerations is the following lemma, which will be useful later.
\begin{lemma}\label{sing}
Let $\nu$ be a Borel measure with compact support in $\R\setminus\{0\}$.  Then the function
\be\la{m28}
f(z)=\intmii \left(\frac{1-z}{1+z}\right)^{is}d\nu(s)
\ee
is in $\dot\HH^{\pul}(D)$\,if and only if $\nu$ is absolutely continuous, with a square-integrable density.
\end{lemma}

\sm
\noindent
{\bf Proof}\\
We can work with the variable $w$ given by~\rf{m10}. Then by our assumptions the functions $f(w)$ given by~\rf{m13n} is in $\dot \HH^{\pul}(\OO)$. This means that $f'(w)$ is in $L^2(\OO)$, and since $f$ is holomorphic in $\OO$, it means that the restriction of $f'$ to the real line is square integrable on the real line. Thus $f'(w_1)$ has the Fourier representation
 \be\la{m30nn}
 f'(w_1)=\int_{-\infty}^{\infty}\psi(s)e^{isw_1}\,ds
 \ee
 for some $\psi\in L^2(\R)$.
 At the same time,
\be\la{m30n}
f'(w_1)=\int_{-\infty}^{\infty}e^{isw_1}is\,d\nu(s)\,,
\ee
and we see that $d\nu(s)=\psi(s)\,ds$ by the Fourier representation uniqueness.
\cara

We note that the functions
\be\la{geneigM}
h(z,\lambda)=\left(\frac{1-z}{1+z}\right)^{\lambda}
\ee
satisfy
$$
\pul(z^2-1)\frac{\partial}{\partial z}\, h(z,\lambda)=\lambda \,h(z,\lambda)
$$
and can be thought of as generalized eigenfunctions of the operator $\pul(z^2-1)\frac{\partial}{\partial z}$ in $\dot \HH^{\pul}$.

\subsection{Generalized eigenfunctions of the operator $L$}\label{ge}
For the spectral analysis of $L$ in the holomorphic sector, we find the analogues of the generalized eigenfunctions~\rf{geneigM} when the simple operator $M$ is replaced by $L$.
The corresponding equation is obtained from~\rf{l19}:
\be\la{g1}
2\lambda z\eta +(z^2-1) (z\eta+zv)'-2z^2(\eta+v)=0\,,\quad \eta=-zv'\,.
\ee
In terms of the Fourier coefficients $\eta_k$ the equation is equivalent to
\be\la{g2}
\begin{array}{rcl}
\lambda \eta_1 & = & A_2\eta_2\\
\lambda \eta_2 & = & A_3\eta_3\\
\lambda \eta_3 & = & B_2\eta_2+A_4\eta_4\\
 &\dots &
\end{array}
\ee
If we choose $\eta_2\ne 0$, then the first equation determines $\eta_1$ (except when $\lambda=0$, of course), and for $\eta_3,\eta_4,\dots$ we get
\be\la{g3}
\eta_{k+1}=A_{k+1}^{-1}\left(\lambda\eta_k-B_{k-1}\eta_{k-1}\right)\,,\qquad k=3,4,\dots
\ee
The usefulness of equation~\rf{g1} is that it enables us to get some control over the functions given by the coefficients calculated from these recursive relations.

As $\eta_0=0$ and $v_0=0$, we can set $\eta=z f$ and $v=z F$. Then
\be\la{g4}
f=-(zF)'
\ee
and~\rf{g1} gives
\be\la{g5}
z(z^2-1)F''+(z^2+2\lambda z -3) F' +2\lambda F=0\,.
\ee
The equation can also be written as
\be\la{g6}
F''+\left[\frac {-1+\lambda}{\,\,\,z-1} + \frac{-1-\lambda}{\,\,\,z+1}+\frac 3 z\right]F'+\frac{2\lambda}{z(z^2-1)} F=0\,.
\ee
This is a classical complex ODE, a special case of the Heun equation~\cite{Heun}. It has four singular points: $z=-1, \,z=0,\,z=1,\,$ and $z=\infty$. All these points are {\it regular singular points}, see for example Chapter X of~\cite{Whittaker} for precise definitions. The local behavior near $z=\infty$ can be investigated, as usual, by setting $z=1/\zeta\,$, which gives
\be\la{g7}
\frac {d^2 F}{d\zeta^2}+\left[\frac{-1+\lambda}{\,\,\,\zeta-1}+\frac{-1-\lambda}{\,\,\,\zeta+1}+\frac 1\zeta\right]\frac{dF}{d\zeta}+\frac {2\lambda}{\zeta(1-\zeta^2)}F=0\,.
\ee
We are interested in solutions which are holomorphic in the neighborhood of $z=0$. To analyze the behavior of the solutions at the other regular singular points, we can use the Frobenius  method and seek solutions in the form
\be\la{g8}
(z-z_0)^r\left(1+a_1(z-z_0)+a_2(z-z_0)^2+\dots\right)\,,
\ee
see for example Chapter X of~\cite{Whittaker} for details.
We will assume that $\lambda$ is a non-zero purely imaginary number. The equation for $r$ (often called the {\it indicial equation}) in our special case is
\be\la{g9}
r(r+\alpha -1)=0\,,
\ee
where $\frac{\alpha}{z-z_0}$ is the term in the square bracket in~\rf{g6} corresponding to the singular point we are considering. In our case $r=0$ is always a solution, and we have one holomorphic solution (up to a multiple, of course) in a neighborhood of each singular point.
The other solution will be crucial for us at the points $-1, 1$. We will also need more information at $z=\infty$.

\sm
\noindent
(i) The general solution near $z_0=1$ when $\lambda$ is not an integer can be expressed as
\be\la{g10}
F=A \,(1-z)^{2-\lambda}U(z,1,\lambda)+B \,V(z,1,\lambda)\,,
\ee
where $A,B\in\CC$,
\be\la{g11}
\begin{array}{rcl}
U(z,1,\lambda) & = & 1+U_1(1,\lambda)(z-1)+U_2(1,\lambda)(z-1)^2+\dots\,\,,\\
 & & \\
V(z,1,\lambda) & = & 1+V_1(1,\lambda)(z-1)+V_2(1,\lambda)(z-1)^2+\dots\,\,,
\end{array}
\ee
where the radius of convergence of the series is at least $R=1$ (the distance between $z_0=1$ and the closest of the remaining singular points, which is $z_0=0$).
The function $z\to(1-z)^{2-\lambda}$ is interpreted in a usual way, along a suitable branch over $\CC\setminus\{1\}$.

\sm
\noindent
(ii) In a similar way, near $z_0={-1}$ we can write
\be\la{g11b}
F=A\,(1+z)^{2+\lambda}U(z,-1,\lambda)+B\,V(z,-1,\lambda)\,.
\ee
Note that the equation has a symmetry $(F(z),\lambda) \to (F(-z),-\lambda)$. In general, the Heun equation has a rich symmetry group, see for example~\cite{Heun}.

\sm
\noindent
(iii) At $z_0=\infty$, which is of course the same as $\zeta=0$ in~\rf{g7}, the indicial equation is $r^2=0$, and there is one holomorphic solution $U$ with $U(0,\lambda)=1$. The general solution is of the form
\be\la{g12}
F=A\,U(\zeta,\lambda)+B\,\left[U(\zeta,\lambda)\log \zeta +V(\zeta,\lambda)\right]\,,
\ee
where $V$ is also holomorphic.

\begin{lemma}\label{noeigenvalues}
The operator $L$ in $X=\dot \HH^{\frac 32}/\CC e_1$ has no point spectrum.
\end{lemma}

\noindent
\sm
We recall that our convention is that constants are factored out already in $\dot\HH^{\frac 32}$.

\sm
\noindent
{\bf Proof}\\
It is easy to see from~\rf{l8} that the kernel of $L$ in $X$ is trivial. (The solutions $\eta=A\sin(\theta-\theta_0)$ are factored out by the projection to $X$.) Therefore we can only consider the case $\lambda=is$ for $s\in\R\setminus\{0\}$. As the functions $(1\pm z)^{is}$ are not in $\dot \HH^{\pul}(D)$ for any $s\ne 0$, the only possibility for the eigenfunctions $\eta$ would be that the corresponding solution $F$ of~\rf{g5} be holomorphic in $\CC$. However, in that case one has to have $B=0$ in the representation~\rf{g12}, which implies that $F$ must be bounded. Hence $F$ is constant by the Liouville theorem, and the claim follows easily.
\cara

\medskip
\noindent
{\bf Remark}\\
Although the operator $L$ does not have any eigenfunctions in $X$, the above analysis shows that the eigenfunctions defined by the formulae~\rf{g3} are regular in $D\setminus\{1,-1\}$, with behavior $\sim (1-z)^{1-is}$ and $(1+z)^{1+is}$ at $1$ and $-1$ respectively.
Such functions just narrowly miss $\HH^{\frac 32}$, and do belong to Sobolev spaces with any lower regularity. The linearized equation therefore has a large set of periodic and almost periodic solutions in spaces just below $\HH^{\frac 32}$. We conjecture that this extends to the non-linear level.

\subsection{Absolute continuity of the spectral measure $\mu$ for the linearized operator $L$.}\label{mains}

Our goal is to prove the following result.
\begin{theorem}\label{mainth}
The measure $\mu$ in the spectral representation~\rf{l20} of $L$ is absolutely continuous.
\end{theorem}

\medskip
\noindent
{\bf Proof}
\\
Let us consider the map $T^{-1}$, where $T$ is the isometry $X\to L^2(\R,\mu)$ defined just before~\rf{l20}. We will represent $\eta=T^{-1}f$ by the Fourier components $\eta_2,\eta_3,\dots$:
\be\la{p1}
f\to T^{-1}f\sim (\eta_2(f),\eta_3(f),\dots)
\ee
For each $k\ge 2$ the map $f\to\eta_k(f)$ is clearly a continuous linear functional on $L^2(\R,\mu)$, and hence we have a representation
\be\la{p2}
\eta_k=\eta_k(f)=\int_{-\infty}^{\infty} f(s) G_k(s)\,d\mu(s)\,.
\ee
with \be\la{p2a}
G_k\in L^2(\R,\mu)\,,\qquad \intmii G_k(s)\bar G_l(s)\,d\mu(s)\sim \frac {\delta_{kl}}{c_k}\qquad \hbox{(no summation)}\,,
\ee
where $c_k$ is given by~\rf{l11}.
Note that
\be\la{p2b}
\eta_k=\frac 1{\,\,c_k}\,(\eta,e_k)_X=\frac 1{\,\,c_k}\,(T\eta,Te_k)_{L^2(\R,\mu)}=\frac 1 {\,\,c_k}\, \intmii f(s) \,\overline{Te_k(s)}\,d\mu(s)\,,
\ee
and hence
\be\la{p2c}
G_k=\frac 1 {\,\,c_k}\, \overline {Te_k}\,.
\ee
Letting  $G=(G_2,G_3,\dots)$, considering $G$ simply as an element of the linear space of infinite sequences of $L^2(\R,\mu)$-functions,
we can write
\be\la{p3}
\eta=\intmii f(s)G(s)\,d\mu(s)\,,\qquad L\eta=\int_{-\infty}^{\infty} f(s)LG(s)\,d\mu(s)\,,
\ee
where the integrals are defined component-by-component.
(Note that $(LG)_k=B_{k-1}G_{k-1}+A_{k+1} G_{k+1}$ is well-defined.)
At the same time, we have by the construction of the isomorphism~$T$,
\be\la{p4}
L\eta=\intmii is f(s) G(s)\,d\mu(s)\,.
\ee
As $f$ was an arbitrary element of $L^2(\R,\mu)$, we see by comparing~\rf{p3} and~\rf{p4} that
\be\la{p5}
L G(s)=is G(s)\,, \qquad \hbox{for $\mu-$almost every $s$.}
\ee
We set
\be\la{p6}
G_1(s)=\frac1{is} A_2G_2(s)\,,\qquad \Psi(z,s)=G_1(s)z+G_2(s)z^2+G_3(s)z^3\dots\,,
\ee
where the definition of $G_1$ reflects the first equation of~\rf{g2}.
Using~\rf{p2a}, we see easily that when $\sigma<1$, the function $z\to\Psi(z,s)$ belongs to $\HH^\sigma(D)$ for $\mu-$almost every $s$. (This is not optimal, but is enough to show that $\Psi(z,s)$ is well-defined as a function of $z$ for $\mu-$almost every $s$.)
Due to~\rf{p5}, the function $z\to\Psi(z,s)$ satisfies equation~\rf{g1} with $\lambda=is$ for $\mu-$almost every~$s$.

Let $\Phi(z,\lambda)$ be a solution of~\rf{g1} defined by the recursive relations~\rf{g3}, with the normalization $\eta_2=1$. By considerations of subsection~\ref{ge}, near $z=1$ we have
\be\la{p7}
\Phi(z,\lambda)=A(\lambda,1)(1-z)^{1-\lambda}+W(z,1,\lambda)\,,
\ee
where $A$ is analytic in $\lambda$ and $W$ is analytic in $z\in D$ and $\lambda$ and $W(\,\cdot\,,1, \lambda)\in \HH^{\frac32}(D)$ as long as $\lambda=is$ and $s\ne 0$. Similarly, near $z=-1$ we have
\be\la{p8}
\Phi(z,\lambda)=A(\lambda,-1)(1+z)^{1+\lambda}+W(z,-1,\lambda)\,,
\ee
with the same properties of $A$ and $W$.
Due to the analyticity, the functions $A(\lambda,\pm 1)$ can vanish only on a countable set of values of the parameter $\lambda$.

As $z\to\Phi(z,is)$ and $z\to\Psi(z,s)$ satisfy the same ODE for $\mu-$almost every $s$ and $\Phi$ is normalized by $\eta_2=1$, we must have
\be\la{p9}
\Psi(z,s)=G_2(s)\Phi(z,is)\,,\quad \hbox{for $\mu-$almost every $s$\,. }
\ee
Assume now that $\mu$ is not absolutely continuous, and let us choose a compact set $E\subset\R\setminus\{0\}$ with vanishing Lebesgue measure such that $\mu(E)>0$ and such that the functions $|G_2(s)|, |A(is,1)|$ and $|A(is,-1)|$ are all bounded below on $E$ by a positive $\ve>0$.
This is possible by the analyticity and non-triviality of the functions $A(is,\pm 1)$, as we already know that $\mu$ contains no Dirac masses, by Lemma~\ref{noeigenvalues}. Let
\be\la{p10}
\mu_E=\chi_E\mu\,,
\ee
where $\chi_E$ is the characteristic function of $E$.

By Lemma~\rf{sing}, for any bounded $\mu-$measurable function $h$ which does not vanish $\mu-$almost everywhere on $E$, the function
\be\la{p11}
F(z)=\intmii \left(\frac{1-z}{1+z}\right)^{-is}\,h(s)\,d\mu_E(s)
\ee
does not belong to $\dot \HH^{\pul}(D)$. Fix such an $h$.
As the function $F$ is smooth away from $z=\pm 1$, the loss of regularity must happen locally at least at one of the points of the set $\{1,-1\}$. Assume it is $z=1$, the other case being essentially the same.

Near $z=1$ we can write
\be\la{p12}
(1+z)^{is}=2^{is}+(1-z)H(z,s)\,,
\ee
where $H$ is analytic in $z$ and $s$ (for $z$ close to $1$).
As the function
\be\la{p13}
z\to (1-z)^{1-is}
\ee
does belong to $\dot \HH^{\pul}(D)$, with a bound on the norm which is uniform in $s\in E$, we conclude that the function
\be\la{p14}
z\to\intmii(1-z)^{-is}h(s)2^{is}\,d\mu_E(s)
\ee
will not belong to $\dot \HH^{\pul}(D)$.
The proof of the theorem is now easily finished by taking
\be\la{p15}
f(s)=\frac{\chi_E(s)h(s)2^{is}}{G_2(s)A(is,1)}
\ee
in formula~\rf{p2}. \cara

\bigskip
\noindent
\begin{corollary}
Let $\eta_0\in \dot \HH^{\frac 32}$ and let $\eta(t)$ be the solution of the linearized equation~\rf{l3} with $\eta(0)=\eta_0$. Then the functions $\eta(t)$ approach zero weakly in $\dot \HH^{\frac 32}/\CC e_1$ as $t\to\pm\infty$.
\end{corollary}

\medskip
\noindent
{\bf Proof}\\
Once we know that $\mu$ is absolutely continuous, the statement follows easily from the Riemann-Lebesgue lemma, the representation~\rf{l20b}, and the fact that the evolution is unitary in $\dot \HH^{\frac 32}/\CC e_1.$
\cara

\bigskip
\noindent
{\bf Remark}
\\
Let $\Phi(z,\lambda)$ be as in~\rf{p7}, and for a smooth $\vf$ compactly supported in $\R\setminus \{0\}$  set
\be\la{p50n}
f(z)=\int_{-\infty}^{\infty}\vf(s)\Phi(z,is)\,ds\,.
\ee
It is natural to expect that
\be\la{p51n}
||f||_{X}^2=\int_{-\infty}^{\infty}|\vf(s)|^2\,\rho(s)\,ds
\ee
for all such $\vf$, with $\rho$ analytic and strictly positive in $\R\setminus\{0\}$.
Our method above can be used to obtain (after more detailed considerations) that we have~\rf{p51n} with $\rho\in L^\infty_{\rm loc}(\R\setminus\{0\})$, and $\rho(s)>0$ for almost every $s$. The density of the functions of the form~\rf{p50n} in $X$ follows
from the proof of Theorem~\ref{mainth}.

\subsection{Operator $L$ in other function spaces}
In this section we will use the similarity of the linearized equation~\rf{l3} with the equation~\rf{x2}. This will allow us to obtain decay estimates for suitable classes of solutions $\eta$ in weighted $L^2$ spaces. When dealing with the linearized operator, it is natural  to work with spaces over $\C$. In the considerations below the functions are complex-valued (unless stated otherwise).

Let us fix $\gamma\in(\frac32, 2)$ and define

\be\la{yzdef}
Y_0=\left\{
f\in L^2(\so)\,,\,
\int_{-\pi}^\pi |f(\theta)|^2|\sin(\theta/2)|^{-2\gamma}<+\infty
\right\}\,,
\ee
On the space $Y_0$ we will take the natural norm
\be\la{nyzdef}
\nyz{f}=\nlt{\,|\sin(\theta/2)|^{-\gamma}\,f}\,.
\ee
 We also define the space
 \be\la{ydef}
 Y=Y_0\oplus \{a+b\sin\theta\,,\,a,b\in\C\}.
 \ee
 with the norm defined for $f\in Y_0$ and $a,b\in \R$ by
 \be\la{nydef}
 ||f+a+b\sin\theta||^2_Y=||f||^2_{Y_0}+|a|^2+|b|^2\,.
 \ee

 For a use of unisotropic Sobolev space for the analysis of the spectral properties of Morse-Smale flows and their action on differential forms (which is in some sense dual to the flow defined by~\rf{M3}) we refer the reader to~\cite{Riviere}.

 Note that each $f\in Y$ has a unique representation of the form
 \be\la{frepres}
 f=f_0+a+b\sin\theta\,\qquad f_0\in Y_0\,,\,\,\,a,b\in\C.
 \ee
 Although values of a function $f\in Y$ may not be defined for all $\theta$, it is natural to define the values of $f$ and $f'$ at $\theta=0$ as
 \be\la{fzvalues}
 f(0)=a\,,\quad f'(0)=b\,.
 \ee
 It is easy to see that this definition agrees with the usual definition when $f\in Y$ is a $C^1$ function. We have the natural projection
 $P_0\colon Y\to Y_0$ defined by for $f_0\in Y_0$ and $a,b\in\R$ by
 \be\la{projyyz}
 P_0(f_0+a+b\sin\theta)=f_0\,.
 \ee
 For $C^1$-functions in $Y$ this amounts to
 \be\la{projco}
 (P_0f)(\theta)=f(\theta)-f(0)-f'(0)\sin\theta\,.
 \ee

 To take advantage of the commutation of $L$ with the Hilbert transform $H$, the following lemma will be useful.
 \begin{lemma}\label{hilbertonY}
 The Hilbert transform $H$ is a continuous operator from $Y$ to $Y$.
 \end{lemma}

\medskip
\noindent
{\bf Proof}

\noindent
Recalling that
\be\la{hp1}
Hf(\theta)=\frac1{2\pi} \int_{-\pi}^{\pi}f(\vartheta)\cot\left( \frac{\theta-\vartheta}2\right)\,d\vartheta\,,
\ee
we see from~\rf{projco} and the fact that $H$ leaves the 2d subspace $\{a\cos\theta+b\sin\theta\}$ invariant that it is enough to show that for a smooth $f\in Y_0$ the function
\be\la{hp2}
P_0Hf:\theta\to \frac 1{2\pi} \int_{-\pi}^{\pi}f(\vartheta)\left[\cot\left(\frac{\theta-\vartheta} 2\right)+\cot\left(\frac{\vartheta}2\right)+\frac {\sin\theta}{2\sin^2\left(\frac{\vartheta}2\right)}\right]\,d\vartheta
\ee
is in $Y_0$, with the corresponding estimate. Using
\be\la{trigid}
\cot\left(\frac{\theta-\vartheta} 2\right)+\cot\left(\frac{\vartheta}2\right)+\frac {\sin\theta}{2\sin^2\left(\frac{\vartheta}2\right)}=\frac{\sin^2\frac \theta2}{\sin^2\frac\vartheta2}\cot\left(\frac{\theta-\vartheta}2\right)\,,
\ee
we see that
\be\la{hp3}
P_0Hf(\theta)=\frac1{2\pi}\int_{-\pi}^{\pi} f(\vartheta)\,\,
\frac{\sin^2\frac \theta2}{\sin^2\frac\vartheta2}
\cot\left(\frac{\theta-\vartheta}2\right)\,d\vartheta\,\,.
\ee

 Let us now write
 \be\la{newcoord}
 f(\theta)=|\sin(\theta/2)|^\gamma g(\theta)\,,\qquad P_0Hf(\theta)=|\sin(\theta/2)|^\gamma G(\theta)\,.
 \ee
 We have to show that the operator
 \be\la{operatorgG}
 g\to G:\theta\to \frac 1{2\pi} \int_{-\pi}^{\pi} g(\vartheta) \frac{|\sin(\theta/2)|^{2-\gamma}}{|\sin(\vartheta/2)|^{2-\gamma}}
 \,\,\cot\left(\frac{\theta-\vartheta}{2}\right)\,d\vartheta
 \ee
 is continuous on $L^2(\so)$.
 This can be either done directly by modifying the proof for the Hilbert transform, or one can use results from the theory of $A_2-$weights. Let us illustrate the latter approach in the case of the Hilbert transform on the real line, leaving the easy adaptation of the proof to $\so$ for the reader. We would like to show that for $\alpha\in[0,1/2)$ the operator $T\colon C^\infty_0\to {\mathcal D}'$ given by
 \be\la{operonr}
 T(x)=\int_{\R} K(x,y)f(y)\,dy\,,\qquad K(x,y)=\frac{|x|^\alpha}{|y|^\alpha}\,\,\frac 1{x-y}
 \ee
 can be continuously extended to an operator from $L^2(\R)$ to $L^2(\R)$.
This is an immediate consequence of the fact that the function $w(x)=|x|^{2\alpha}$ is an $A_2-$weight, see Stein~\cite{Stein} for the definitions. To verify this, we note that the function $f$ belongs to $L^2(\R)$ if and only if the function $y\to f(y)|y|^{-\alpha}$ belongs to $L^2(\R, w(y)\,dy)$ and, similarly, a function $x\to |x|^\alpha g(x)$ belongs to $L^2(\R)$ if and only if $g\in L^2(\R, w(x)\,dx)$.
\cara

\bigskip
\noindent
Above we complemented the space $Y_0$ by the 2-dimensional space of the functions of the form $a+b\sin\theta$. When working with holomorphic functions, it is better to work with another natural complementary space of $Y_0$, defined as
\be\la{Zdef}
\Z_0=\{f\,,\,f(z)=a+b(z-1)\,,\,a,b\in\C\}.
\ee
Clearly $Y_0\cap \Z_0=\{0\}$ and $Y=Y_0\oplus \Z_0$. Moreover, the functions in $\Z_0$ are holomorphic. for $g\in Y_0$ and $f\in Y$ given by $f(z)=g(z)+a+b(z-1)$ we can define a norm in $Y$ (equivalent to the previously defined norm $||\,\cdot\,||_{Y}$) by
\be\la{Ynorm}
||f||^2=||g||^2_{Y_0} + |a|^2 + |b|^2\,.
\ee
In what follows we will not distinguish too carefully between the two equivalent norms, as the distinction is not important for our purposes.
For $f(z)=g(z)+a+b(z-1)$ with $g\in Y_0$ we can define (keeping in mind that $z=e^{i\theta}$)
\be\la{fonedef}
f(1)=a\,,\qquad \frac {\partial f}{\partial z}(1)=b\,,
\ee
where the value $1$ of course refers to the variable $z$, and we use $iz\frac{\partial}{\partial z}=\frac{\partial}{\partial \theta}$.
The decomposition $Y=Y_0\oplus \Z_0$ defines projections
\be\la{projPQ}
P\colon Y\to Y_0\,,\qquad Q\colon Y\to \Z_0\,.
\ee
It is perhaps worth emphasizing that $P$ does not coincide with $P_0$ defined above.

The proof above of the continuity of the Hilbert transform on $Y$ can we re-written in the complex notation, using
\be\la{complexcot}
\cot\left(\frac{\theta-\vartheta}2\right)=i\,\frac{z+w}{z-w}\,,\qquad z=e^{i\theta}\,,\,w=e^{i\vartheta}
\ee
and replacing~\rf{trigid}  by
\be\la{trigid2}
\frac{z+w}{z-w}
-\frac{1+w}{1-w}
+\frac{2w\,(z-1)}{(1-w)^2}
=
\frac{2w\,(1-z)^2}{(z-w)(1-w)^2}\,.
\ee

\bigskip
Equation~\rf{l19} for $\eta$ can be rewritten as
\be\la{etaeqz}
\eta_t+\pul(z^2-1)\eta'-z\eta -\pul (z+\frac 1z)v=0\,,\quad \eta=-zv'\,,\quad v|_{z=0}=0\,.
\ee
Note that when considering $\eta$ as a function on $\so$, the equation makes sense even when $\eta$ is not holomorphic, as for functions on $\so$ we have $iz\frac{\partial}{\partial z}=\frac{\partial}{\partial\theta}$. However, the equation $\eta=-zv'$ coincides with $v_\theta=H\eta$ only on holomorphic functions.
We first analyze the equation
\be\la{simpleeq}
\eta_t+\pul(z^2-1)\eta'-z\eta=0,
\ee
which is equivalent (on the circle $\so$) to $\xi_t+\sin\theta\,\,\xi_\theta-\cos\theta\,\,\xi=0$ via the change of variables $\eta=z\xi$.
The equation can be solved explicitly as follows. The flow map generated by the ODE
\be\la{ode1}
\dot z=\pul(z^2-1)
\ee
is
\be\la{sol1}
\phi_t(z)=\frac{z-\tau}{1-\tau z}\,,\qquad \tau=\tanh\frac t2\,.
\ee
Hence the solution of~\rf{simpleeq} with the initial condition $\eta|_{t=0}=\eta_0$ is
\be\la{sol2}
\eta(z,t)=\phi'_t(\phi^{-1}_t(z))\,\,\eta_0(\phi^{-1}_t(z))
\ee
and $\nyz{\eta(t)}^2$ is given by
\be\la{decay1}
\nyz{\eta(t)}^2=\iso|\eta(z,t)|^2|z-1|^{-2\gamma}\,d\mathscr H^1(z)\,,
\ee
where $\mathscr H^1$ is the standard 1d measure on the circle.
Setting $z=\phi_t(w)$ in the last integral and using~\rf{sol2} together with
\be\la{subst1}
d\mathscr H^1(z)= |\phi_t'(w)|\,d\mathscr H^1(w)=\frac{1-\tau^2}{|1-\tau w|^2}\,d\mathscr H^1(w)\,,
\ee
and
\be\la{pomocna}
|w-\tau|=|1-\tau w|\ge |1-\tau|\,\quad\hbox{when $ |w|=1$ and $0\le\tau<1$}\,
\ee
we see that
\be\la{calc1}
\nyz{\eta(t)}^2\le(1-\tau)^{2\gamma-3}\iso |\eta_0(w)|^2|w-1|^{-2\gamma}
d\mathscr H^1(w)
\le e^{-2\beta t}\nyz{\eta_0}^2\,,\quad \beta_0=\gamma-\frac32\,.
\ee
Let us denote
\be\la{lz}
L_0\eta=-\pul(z^2-1)\frac{\partial\eta}{\partial z} +z\eta
\ee
We can state~\rf{calc1} as follows.
\begin{lemma}\label{decaylz}
\be\la{decayest}
\nyz{e^{tL_0}\eta_0}\le e^{-\beta_0 t}\nyz{\eta_0}\,\quad t\ge 0\,.
\ee
\end{lemma}

\medskip
\noindent
{\bf Proof}
\,\, See above.\,\,\cara

\bigskip
Our goal is to prove a suitable version of this estimate for the operator $L$.

Let us define an operator $\eta\to K\eta$ on $L^2(\so)$ by
\be\la{k}
K\eta=
(\cos\theta)\,v\,,\quad v_\theta=H\eta\,,\quad \int_{-\pi}^{\pi} v(\theta)\,d\theta=0\,.
\ee
Note that on holomorphic functions
\be\la{khol}
K\eta=\pul(z+\frac 1z)v\,,\quad \eta=-zv'\,,\quad v|_{z=0}=0\,.
\ee
It is worth noting that
\be\la{koz}
K\cdot 1=0\,,\quad K\cdot z= -\pul(z^2+1)\,.
\ee
\begin{lemma}\label{Kcompact}
The operator $K$ maps $Y$ to $Y$ and is continuous and compact as an operator on~$Y$.
\end{lemma}
{\bf Proof}\,\,\,
The statement is clearly true for the restriction of $K$ to $\Z_0$. Therefore is is enough to show that $K\colon Y_0\to Y$ is continuous and compact. As the maps $f\to (\cos\theta) f$ and $f\to Hf$ are continuous from $Y$ to $Y$, we only have to show that taking primitive of a function in $Y_0$ with zero average is a compact map from $Y_0$ to $Y$. Let us write elements of $Y_0$ as $f(\theta)=|\sin\frac\theta2|^\gamma\,g(\theta)$ with $g\in L^2(\so)$.
Letting $w(\theta)=|\sin\frac\theta2|^\gamma$, define
\be\la{prim}
Tg(\theta)=\int_0^\theta \frac{w(\vartheta)}{w(\theta)}g(\vartheta)\,d\vartheta\,.
\ee
We need to show that $T$ is a compact operator from the subspace of $L^2$ of functions with $\int_{-\pi}^\pi w(\theta)g(\theta)=0$ into $L^2$.
Let us fix some $0<\theta_0<\frac\pi2$. For $|\theta|\le\theta_0$ we have
\be\la{cauchyschwartz}
|Tg(\theta)|^2
\le||g||_{L^2}^2
\int_0^\theta\frac{w^2(\vartheta)}
{w^2(\theta)}\,d\vartheta\,.
\ee
This gives
\be\la{compest}
|Tg(\theta)|^2\less ||g||_{L^2}^2|\theta|\,,\qquad |\theta|\le \theta_0\,.
\ee
This gives sufficient control near $\theta=0$. In regions  away from small neighborhoods of $\theta=0$  we can use standard results about compactness of integral operators.
\cara

Let $Z_0$ be the subspace of $Y$ spanned by $1,\cos\theta$, and $\sin\theta$. Recalling that $K\cdot 1=0$, and the equilibria~\rf{e1}, we see that
\be\la{kzinv}
L\cdot 1=\cos\theta\,,\,\quad L\cdot \cos\theta=0\,,\quad L\cdot \sin \theta = 0\,.
\ee
(Of course, the last two equalities can be also seen by a direct substitution of $\cos\theta$ and $\sin\theta$ into $L$.) Hence $Z_0$ is invariant under $L$ and $L$ is well-defined on the factor space
\be\la{factorspace}
\tY=Y/Z_0\,.
\ee
\begin{lemma}\la{maindecay}
For each $\beta<\beta_0=\gamma-\frac 32$  we have the decay estimate
\be\la{maindecayest}
||e^{tL}\eta||_{\tY} \le Ce^{-\beta t}||\eta||_{\tY}\,,\quad \eta\in Y\,,
\ee
where $C=C(\beta)$ is a suitable constant.
\end{lemma}

This implies that for $|k|\ge 2$ the Fourier coefficients $\eta_k(t)$ of $e^{tL}\eta$ decay exponentially as $t\to\infty$. This is a considerable strengthening of our analysis in $\dot H^{\frac32}$ in the previous subsection. It is this exponential decay (together with the $\dot H^{\frac 32}$ estimate), which will enable us to do perturbation analysis near equilibria in the non-linear case.

\medskip
\noindent
{\bf Proof of the lemma}\,\,
We note that the commutator $[L,H]$ (where $H$ is the Hilbert transform) vanishes on $\tY$ (although it does not vanish on $Y$, as we now do not assume that $\iso\eta =0$), and hence we can decompose $\eta$ as in~\rf{l14} and prove the decay separately for $\eta_+$ and $\eta_-$. For the rest of the proof we will assume that $\eta=\eta_+$ is holomorphic in the unit disc. The subspace of holomorphic functions in $Y$ will be denoted by $\YY$. The factor space $\YY/\Z_0$ will be denoted by $\tYY$. (Recall that $\Z_0$ is the linear space of $1$ and $z$.) The projections $P$ and $Q$ defined by the decomposition $Y=Y_0\oplus\Z_0$ (see~\rf{projPQ}) map $\YY$ into itself, and we will denote their restrictions to $\YY$ also by $P$ and $Q$, respectively. We will also denote by $\YY_0$ the holomorphic functions in $Y_0$. Clearly $P\YY=\YY_0$.

Letting
\be\la{LoKodef}
L_1=L_0P\,,\quad K_1=L_0Q+K
\ee
we clearly have
\be\la{Ldecomp}
L=L_1+K_1\,.
\ee
Moreover, as $K$ is compact by Lemma~\rf{Kcompact} and $Q$ has a finite-dimensional range on which $L_0$ is continuous, the operator $K_1$ is compact in $Y$. In view of Lemma~\rf{decaylz}, we also have
\be\la{decayL0}
||e^{tL_1}\eta||_{\tYY}\le Ce^{-\beta t}||\eta||_{\tYY}\,,\qquad \eta\in \YY\,.
\ee
In this situation the only obstacle to the decay  estimate~\rf{maindecayest} can come from possible points of the spectrum of $L$ in the region $\{\lambda\,,\,\Re\lambda>-\beta_0\}$ see~\cite{Engel-Nagel}, Section 2 of Chapter IV, Corollary 2.11 and Proposition 2.12. Hence we need to study the solutions $\eta\in\YY$ of $L\eta=\lambda\eta$. Using this equation, we can again look at the recursive relations~\rf{g3} for the Fourier coefficients and conclude that $\eta$ has to be a holomorphic function satisfying~\rf{g1}.
We return to the analysis of the equation~\rf{g1} for the holomorphic eigenfunctions in subsection~\ref{ge}, which we re-write here for the convenience of the reader:
\be\la{rE}
z(z^2-1)F''+(z^2+2\lambda z-3)F'+2\lambda F=0\,.
\ee
By the method of Frobenius discussed in subsection~\ref{ge}, for each $\lambda\in\CC$ there is a unique solution of~\rf{rE} holomorhic in the unit disc with $F(0)=1$. We will denote this solution $F(z,\lambda)$.

We will use the following notation: if $g(z)=\sum_{k=0}^\infty g_kz^k$ is a function which is holomorphic in a neighborhood of $z=0$, we define $\bar g(z)=\sum_{k=1}^\infty \bar g_k z^k$. Clearly $\bar g$ is holomorphic, with the same radius of convergence for its Taylor series at $z=0$ and $\overline{g(z)}=\bar g(\bar z)$. Applying this notation to the function $z\to F(z,\lambda)$ where $\lambda$ is considered as a parameter, we can write
\be\la{np1}
 \bar F(z,\lambda)=F(z,\bar\lambda)\,.
\ee
and
\be\la{np1b}
\overline {F(z,\lambda)}=F(\bar z, \bar \lambda)\,.
\ee
We also note that $F(-z,-\lambda)$ satisfies~\rf{rE} and its value at $z=0$ is $1$, which means by uniqueness that
\be\la{np2}
F(-z,-\lambda)=F(z,\lambda)\,.
\ee
Let us now assume that $\lambda=is$ for $s\in \R\setminus\{0\}$, and let us look at the function $x\to F(x,\lambda)$ for real $x$. Using~\rf{np1b} and~\rf{np2}, we see that
\be\la{np3}
\overline {F(-x,is)}=F(-x,\overline{is})=F(-x,-is)=F(x,is)\,.
\ee
Given the local form of $F$ near $z=-1$ and $z=1$, we see that for a non-zero purely imaginary $\lambda$, the function $z\to F(z,\lambda)$ is singular at $z=1$ if and only if it is singular at $z=-1$. From this it is easy to see that the generalized eigenfunction $\eta$ corresponding to $\lambda = is$ cannot belong to $Y$, unless it is regular in $\C$. However, we have seen in the proof of Lemma~\ref{noeigenvalues} that such functions have to be constant, which correspond to $\eta(z)=cz$, which project to $0$ in $\tYY$. Hence the imaginary axis does not contain any points of the spectrum of $L$ considered as an operator from $\tilde \YY$ to $\tilde \YY$.

It remains to deal with the case when $\lambda$ is not on the imaginary axis, with $\Re \lambda>-\beta_0$.
  Assume that $\lambda$ is an eigenvalue of $L$ in $\tYY$ with $\Re \lambda>-\beta_0$ and $\Re\lambda\ne0$. The operator $\lambda-L$ is a compact perturbation of the invertible operator $\lambda-L_0P$ (with all operators being considered on $\tYY$), and hence it is Fredholm. The projection $P_\lambda$ on the eigenspace of $\lambda$ (or, equivalently, the kernel of $\lambda-L$) is given by
\be\la{projplambda}
P_\lambda\eta=\frac1{2\pi i} \int_C (z-L)^{-1}\eta\,dz\,,
\ee
where $C$ is a sufficiently small circle around $\lambda$ (taken with the positive orientation). We take $C$ so that it does not intersect the imaginary axis $\R i$\,. When $\eta\in \dot H^{\frac 32}$, the integral vanishes, as the region surrounded by the contour $C$ does not contain any spectral value of $L$ (considered in $\dot H^{\frac 32}$), and hence the integral has to vanish. As smooth functions are dense in $Y$, the integral has to vanish in $\tY$ for any $\eta\in Y$.\cara

\medskip
\noindent
\vbox{\noindent{\bf Remarks}

\smallskip
\noindent
1. The equation~\rf{rE} very likely has non-trivial solutions $F(z,\lambda)$ which are regular both at} $z=0$ and $z=1$ for a countable set of real $\lambda_n\searrow-\infty$ which satisfy $\lambda_n<-\beta_0$. These solutions do not interfere with our estimate.  \\
2. In the definition of the space $Y_0$ we work with approximation of functions by affine maps near a point. There is a natural generalization of $Y_0$ to higher-order approximations. Let $m\ge 1$ be an integer and let $\gamma\in (m-\pul, m)$. Let
\be\la{highersp1}
Y^{(m)}_0=\{f\in L^2(\so)\,;\,|z-1|^{-\gamma}\,f\in L^2(\so)\}\,,
\ee
with the norm $||f||_{\ymz}=||\,|z-1|^{-\gamma}\, f\,||_{L^2}$. We denote by $\PP_{m}$ the set of polynomials of degree $m$ and define
$Y^{(m)}=\ymz\oplus \PP_{m-1}$, with a norm of $f=g+p\,,\,\,g\in\ymz,p\in\PP_{m-1}$  given by $||f||_{\ym}^2=||g||_{\ymz}^2+||p||_{\PP_m}^2$, where $||\,\cdot\,||_{\PP_m}$ is a norm on $\PP_m$. (This does not define the norm uniquely, of course, but any two norms defined in such a way are equivalent.)
Replacing~\rf{trigid2} by its $m-$th order generalization
\be\la{highersp2}
\frac{z+w}{z-w}=\frac{1+w}{1-w}+
\sum_{k=1}^{m-1}\left(\frac{z-1}{w-1}\right)^k\frac{2w}{w-1}+
\left(\frac {z-1}{w-1}\right)^m\frac{2w}{z-w}\,,
\ee
one can use similar arguments as above to show that $\ym$ is invariant in $H$, for $f\in \ym$ the map $f\to  (f(1),f'(1),\dots,f^{(m-1)}(1))$ is continuous, and the evolution by $L_0$ and $L$ preserves $\ym$.

\medskip
\noindent
3. It is not hard to see that the evolution given by $\dot\eta=L\eta$ preserves smooth functions and many other regularity classes, such the above defined $\ym$ spaces. This can be seen in many ways. For example, it is obvious from the explicit formulae that the operator $e^{tL_0}$ preserves various regularity classes. The generator of $L$ is a bounded (and even compact) perturbation of the generator of $L_0$ in the spaces of holomorphic functions discussed above, and this can be used to show the desired regularity, as long as $H$ and $K$ preserve the regularity classes. A useful corollary of this is that when $\eta_0\in Y$, the equations for $\eta(z_0,t)$ and $\eta'(z_0,t)$ are well defined, and are the same as in the case of smooth functions. To see this, one can use the density of the smooth functions in $Y$, continuity of $e^{tL}$ in $Y$, and the continuity of $\eta\to (\eta(1),\eta'(1))$ on $Y$.
The same will be true for the linearization of~\rf{i9j} about an equilibrium.

\bigskip
It will be useful to look at the decay of a general solution of the form~\rf{i9j} of the linearized equation, with the calibration $\tilde u(\theta_0,t)=0$.
Let us set $\theta_0=0$. The equation then is
\be\la{moddec1}
\eta_t+\sin\theta \,\,\eta_\theta-\cos\theta \,\,\eta+\sin\theta \,\,v_\theta- \cos\theta\,\,v +\cos\theta\,\,v(0,t)\,=\,0\,,\quad v_\theta=H\eta\,,\quad \int_{-\pi}^{\pi} v(\theta)\,d\theta=0\,,
\ee
and we assume $\int_{-\pi}^{\pi}\eta(\theta)\,d\theta=0$.
By the previous remark, for any solution with the initial condition $\eta_0\in Y$ (with zero average) the solution will be in $C([0,\infty), Y)$, the values $\eta(0,t),\eta_\theta(0,t)$ will be well defined for all $t$ and the same equations for $\eta(0,t), \eta_\theta(0,t)$ are satisfied as in the case when $\eta$ is smooth. In particular, it is easy to check that the condition $\eta(0,t)=0$ is preserved under the evolution by~\rf{moddec1}, and so is the condition \hbox{$(\eta(0,t),\eta_\theta(0,t))=(0,0)$.} Note that the term with $v(0,t)$ in~\rf{moddec1} does not affect the projection of $\eta$ to $Y/Z_0$, and hence Lemma~\rf{maindecay} implies that for a suitable

\be\la{xiformula}
\xi(\theta,t)=a_1(t)e^{i\theta}+a_{-1}(t)e^{-i\theta}
\ee
we have
\be\la{xidecay}
||\eta(t)+\xi(t)||_{Y}\le Ce^{-\beta t}\,.
\ee
Denoting by $\eta_k(t)$ the Fourier coefficients of $\eta(\theta,t)$, we have
\be\la{omoe}
\begin{array}{ccccccc}
\dot\eta_1\,\,\, & = & \pul(\eta_1+\eta_{-1})& + & A_2\eta_2& + & \pul\sum_{|k|\ge 2}\frac{\eta_k}{|k|}\,\,\\
\dot \eta_{-1} & = & \pul (\eta_1+\eta_{-1})& +& B_{-2}\eta_{-2}& +&
\pul\sum_{|k|\ge 2}\frac{\eta_k}{|k|}\,.
\end{array}\,
\ee
From~\rf{xidecay} we see that the terms not containing $\eta_{1}$ and $\eta_{-1}$ decay as $e^{-\beta t}$ as $\to\infty$.
Hence for $y_1(t)=\eta_1(t)+\eta_{-1}(t)$ we have
\be\la{yeq1}
\dot y_1 = y_1+ g_1(t)\,,\quad |g_1(t)|=O(e^{-\beta t})\,,\,\,\,t\to\infty.
\ee
Therefore
\be\la{ysol1}
y_1(t)=-\int_t^\infty e^{t-s}g_1(s)\,ds + \tilde c_1e^t\,
\ee
for some $\tilde c_1\in \C$. (For real-valued solutions $\eta$ we will have $\tilde c_1\in \R$.)
Similarly, for $y_2(t)=\eta_1(t)-\eta_{-1}(t)$ we gave an equation
\be\la{yeq2}
\dot y_2 = g_2(t)\,,\quad |g_2(t)|=O(e^{-\beta t})\,,\,\,\,t\to\infty,
\ee
with the general solution
\be\la{ysol2}
y_2(t)=-\int_{t}^\infty g_2(s)\,ds + \tilde c_2\,,
\ee
where $\tilde c_2\in \C$ (and $\tilde c_2\in i\R$ for real-valued solutions). We conclude that in~\rf{xiformula} we must have
\be\la{coeftomo}
a_1(t)=c_1e^t+c_2+O(e^{-\beta t})\,,\quad a_{-1}(t)=c_1e^t-c_2+O(e^{-\beta t})\,
\ee
where $t\to\infty$, where $c_1,c_2\in\R$ are suitable constants. If $\eta(0,t)=0$, then $\eta(0,t)=0$ for all $t$ and $c_1=0$. If $\eta_\theta(0,t)=0$, then $\eta_\theta(\theta,t)=0$ for all $t$ and $c_2=0$.
In particular, we have proved the following statement.
\begin{lemma}\label{maindecay2}
Let us denote by $e^{t\tilde L}$ the semigroup in $Y$ generated by equation~\rf{moddec1}. There exists $C_0\in\R$ such that for any $\eta_0\in Y_0$ we have
\be\la{maindecay2f}
||e^{t\tilde L}\eta_0||_{Y_0}\le C_0e^{-\beta t}||\eta_0||_{Y_0}\,.
\ee
\end{lemma}

\section{Nonlinear stability}\label{nlsection}
We will consider the non-linear stability of the steady state $\Om=-A\sin(\theta-\theta_0)$ of equation~\rf{3}. We assume the initial data is of form $\om_0=\Om+\ve\eta_0$, where $\eta_0$ is a sufficiently regular function (roughly of size of order one in a suitable norm) and $\ve$ is small. By a suitable rotation we can assume $\om_0(0)=0$, and we can adjust $\Om$ by changing $A$, in necessary, so that $\Om_\theta(0)=\om_{0\theta}(0)$, which then gives $\eta_0\in Y_0$. After these transformation we can also multiply $\Om$ by a suitable factor and rescale time, so we are in the situation with $\Om=-\sin\theta$ and $\eta\in Y_0$. Our main assumption now will be that
\be\la{nstab1}
||\eta_0||_{\dot H^{\frac 32}}\less 1\,,\quad ||\eta_0||_{Y_0}\less 1\,.
\ee
We will consider the evolution in the gauge~\rf{i9j} with $\theta_0=0$, so that the condition $\eta(0,t)=0$ is preserved during the evolution.
The evolution is given by
\be\la{nstab2}
\eta_t+[U,\eta+\tilde v] + \ve[\tilde v, \eta]=0\,,\quad \tilde v_\theta=H\eta\,,\quad \tilde v(0,t)=0\,.
\ee
Unless otherwise stated, the functions $\eta, v$ etc. in this section are considered to be real-valued.

The non-linear problem is well-posed locally in time and the regularity of the initial data is preserved in the closed time interval $[0,T]$ under the assumptions above, as long as the quantity $\int_0^T||\eta(t)||_{\dot H^{\frac12}}\,dt$ is finite. This is similar to the Beale-Kato-Majda-type criterion for 3d incompressible Euler, see \cite{BKM}, and can be proved along similar lines. A slightly different form, in which the last integral is replaced by $\int_0^T||v_\theta(t)||_{L^\infty}\,dt$,  is proved in~\cite{OSW}. That form is fully sufficient for our purposes here.

For the remainder of this section with will assume $\Om=-\sin\theta\,,\,U=\sin\theta,$ and $\theta_0=0$. As above, we will use the notation $\tL$ for the operator $\eta\to-[U,\eta+\tilde v]$.

Our goal is to prove the following result.
\begin{theorem}\label{main}
When $\ve>0$ is sufficiently small, the local solution $\eta$ of equation~\rf{nstab2} can be continued globally for any initial condition satisfying~\rf{nstab1}, and satisfies
\be\la{global}
||\eta(t)||_{Y_0}\le 2C_0e^{-\beta t}\,\quad\hbox{and $\quad||\eta(t)||_{H^\frac32}\le 2C_0$\,,}
\ee
for all $t\ge 0$, where $C_0$ is the constant from Lemma~\ref{maindecay}.
\end{theorem}

For the proof of the theorem we first establish a few auxiliary results.

\subsection{Energy estimate}\label{energyestimate}
Let us denote by $M$ the operator on functions $\so$ given by the Fourier multiplier \\$\sqrt{(k^2-1)(|k|+1)}$, i.\ e.\
\be\la{nstab3}
\widehat{Mf}_k=\sqrt{(k^2-1)(|k|+1)}\,\,\widehat f_k\,.
\ee
We saw in the previous section (see~\rf{ckprecise}) that the quadratic form $\iso |M\eta(\theta)|^2\,d\theta$ is preserved by $e^{t\tL}$.
The following lemma provides an estimate of how the conservation is affected by the non-linear term.

Recall that $Z_0$ is the subspace of functions on $\so$  generated by  $1,\cos\theta,\sin\theta$, and that for $\eta\in H^s(\so)$ we define $\tilde v$ by the ``Biot-Savart law"~\rf{i9f}.

\begin{lemma}\label{energyest} Given  $\sigma>\pul$, there exists $C=C(\sigma)\ge0$ such that for each $\eta\in H^{2}(\so)$ we have
\be\la{mainenest}
\left|\iso (M[\tilde v,\eta])M\eta\,d\theta\right|\le C||\eta||_{\dot H^{\sigma}/Z_0}||M\eta||^2_{L^2}\,.
\ee
\end{lemma}

\bigskip
\noindent
{\bf Proof} \,\,\, It is not hard to check by direct calculation that the integral on the left-hand side of~\rf{mainenest} does not change when we change $\eta$ by a function in $Z_0$. It is therefore enough to prove the estimate with the right-hand side replaced by $C||\eta||_{H^\sigma}||\eta||^2_{\htt}$.
Recalling that $[\tilde v,\eta]=\tilde v\eta_\theta-\tilde v_\theta\eta$ we note that the part of the left-hand-side of~\rf{mainenest} arising from the term $\tilde v_\theta\eta$ is easily estimated as required by the standard estimate for multiplication of functions in Sobolev spaces
\be\la{sobolev1}
||\tilde v_\theta \eta||_{\htt}\less ||\tilde v_\theta||_{L^\infty}||\eta||_{\htt}+||\tilde v_\theta||_{\htt}||\eta||_{L^\infty}\,.
\ee
As is usual in similar situations, the main point in estimating the part arising from $\tilde v\eta_\theta$ is to integrate by parts. We write
\be\la{comm1}
M\left(\tilde v\eta_\theta\right)=\tilde v \,M\eta_\theta+[M,\tilde v]\,\eta_\theta\,,
\ee
where we denote by $[M,\tilde v]$ the commutator of the operator $M$ and the multiplication operator by $\tilde v$. The term coming from $\tilde v \,M\eta_\theta $ is estimated from
$$
\iso \tv (M\eta)_\theta(M\eta)\,d\theta=-\pul\iso\tv_\theta|M\eta|^2\,d\theta\,,
$$
and it only remains to estimate the term with the commutator. The standard Kato-Ponce estimate (see~\cite{KP}) applied to our situation gives
\be\la{KatoPonce}
||[M,\tv]\eta_\theta||_{L^2}
\less
||\tv_{\theta}||_{L^\infty}
||\eta_\theta||_{H^\pul}+
||\tv||_{\htt}
||\eta_\theta||_{L^\infty}\,,
\ee
which is close to what we need, except for the term $||\eta_\theta||_{L^\infty}$. We can replace $||\eta_\theta||_{L^\infty}$ by $||\eta||_{H^\pul}$ at the cost of adding an extra $\ve-$-derivative to $||\tv||_{\htt}\sim ||\eta||_{H^\pul}$, obtaining
\be\la{maincommutator}
||[M,\tv]\eta_\theta||_{L^2}
\less
||\tv_{\theta}||_{L^\infty}
||\eta_\theta||_{H^\pul}+
||\eta||_{H^\sigma}
||\eta_\theta||_{H^\pul}\,\less\,||\eta||_{H^\sigma}
||\eta_\theta||_{H^\pul}\,,
\ee
which, together with the other estimates above, gives~\rf{mainenest}.

Estimate~\rf{maincommutator} does not have the optimal scaling (unlike the Kato-Ponce estimate) and can be proved by a standard application of Cauchy-Schwatz inequality on the Fourier side. For the convenience of the reader we outline the proof.
Let us write $f=\tilde v\,,\,g=\eta_\theta$, and let $f_k, g_l$ denote the Fourier coefficients of $f$ and $g$, respectively. The Fourier series for the commutator is given (up to a factor of $(2\pi)^2$) by
\be\la{commfourier}
\sum_{k,l} K(k,l)F_kG_le^{i(k+l)\theta}\,,
\ee
where
\be\la{expressions1}
F_k=\langle k\rangle^{1+\sigma} f_k\,,\quad G_l=\langle l\rangle^{\frac 12} g_l\,,\quad K(k,l)=\frac{M(k+l)-M(l)}{\langle k\rangle^{1+\sigma}\langle l\rangle^{\frac12}}\,,
\ee
with the usual notation $\langle k\rangle=\sqrt{1+k^2}$.
Hence, by the Cauchy-Schwartz inequality
\be\la{expressions1a}
\begin{split}
||[M,\tilde v]\eta_\theta||^2_{L^2} & \lesssim \sum_m\left(\sum_{k+l=m}|K(k,l)|^2\right)\left(\sum_{k+l=m}|F_k|^2|G_l|^2\right)\\
& \le \left(\sup_m\sum_{k+l=m}|K(k,l)|^2\right)\sum_{k,l}|F_k|^2|G_l|^2 \\
& = \left(\sup_m\sum_{k+l=m}|K(k,l)|^2\right)\left(\sum_k|F_k|^2\right)\left(\sum_l|G_l|^2\right)\,.
\end{split}
\ee

Letting
\be\la{expression3}
C_*=\sup_m \left(\sum_{k+l=m} |K(k,l)|^2\right)^\frac12\,,
\ee
we see from~\rf{expressions1a} that
\be\la{expression4}
C_* ||f||_{H^{1+\sigma}}||g||_{H^{\frac12}}\,.
\ee
The proof is completed by showing that $C_*$ is finite, which is an easy exercise. The main point is that for large $m$ and $k+l=m$  one has to use the cancellation in $M(k+l)-M(l)$ when $k$ is small relative to $l$.
\cara

\bigskip
Let us now assume that on a time interval $[0,T]$  we have a solution of~\rf{nstab2} with initial condition $\eta_0$ satisfying~\rf{nstab1} such that
\be\la{globalGamma}
||\eta(t)||_{Y_0}\le \Gamma e^{-\beta t}\,\quad\hbox{and $\quad ||M\eta(t)||_{L^2}\le \Gamma\,,\quad t\in[0,T]$\,,}
\ee
where $\Gamma$ is a definite constant.
Then for $\alpha$ slightly below $\frac23$ we have, based on the equation and Lemma~\rf{energyest},
\be\la{estM1}
\begin{split}
\frac {d}{dt} \iso ||M\eta(t)||^2 & =-\iso 2\ve \left(M[\tv,\eta]\right)M\eta\,\,d\theta
\less \ve ||\eta(t)||_{\dot H^\sigma/Z_0}||M\eta(t)||_{L^2}^2\\
& \less \ve\left(||M\eta(t)||_{L^2}^{1-\alpha}||\eta(t)||_{Y_0}^\alpha+
||\eta||_{Y_0}\right)||M\eta(t)||_{L^2}^2\\
&\less \ve\Gamma e^{-\beta\alpha t}||M\eta(t)||_{L^2}^2\,.
\end{split}
\ee
Using $||M\eta(0)||_{L^2}\le 1$ and the Gronwall inequality, we see that
\be\la{estM2}
||M\eta(t)||_{L^2}\le e^{c\ve\Gamma}\,,
\ee
where $c$ is some fixed constant obtained from the various constants involved in the inequalities we have used.

\subsection{Estimates in $Y_0$}
Let us set
\be\la{ye0}
b(t)=\tv_\theta(0,t)=v_\theta(0,t)\,,\quad w=\tv-b(t)U
\ee
and let us write the non-linear term as
\be\la{ye1}
[\tv,\eta]=[w,\eta]+b(t)[U,\eta]\,.
\ee
To estimate $[w,\eta]$ in $Y_0$, we estimate separately $w\eta_\theta$ and $w_\theta\eta$.
\bigskip
For the last term we clearly have
\be\la{ye2}
\begin{split}
||w_\theta \eta||_{Y_0} & \less ||w_\theta||_{L^\infty}||\eta||_{Y_0}=||H\eta-H\eta(0)U_\theta||_{L^\infty}||\eta||_{Y_0}\\
& \less\left(||H\eta||_{L^\infty}+|H\eta(0)|\right)||\eta||_{Y_0}\less
\left(||\eta||_{\htt}^{1-\alpha}||\eta||_{L^2}^{\alpha}+||\eta||_{Y_0}\right)||\eta||_{Y_0}\\
&\less\left(||M\eta||_{L^2}^{1-\alpha}||\eta||_{L^2}^{\alpha}+||\eta||_{L^2}+||\eta||_{Y_0}\right)||\eta||_{Y_0}\\
&\less ||M\eta||^{1-\alpha}||\eta||_{Y_0}^{1+\alpha}+||\eta||_{Y_0}^2\,.
\end{split}
\ee
for an $\alpha$ slightly below $\frac23$. To estimate $w\eta_\theta$, we note that from the proof of Lemma~\ref{Kcompact} and the fact that $w(0)=w'(0)=0$ it follows that
\be\la{ye2b}
|w(\theta)|\less \left(\sin\frac\theta2\right)^2||\eta||_{Y_0},
\ee and therefore
\be\la{ye3}
||w\eta_\theta||_{Y_0}\less||\eta||_{Y_0}||\eta_\theta||_{L^2}\less ||\eta||_{Y_0}||\eta||_{L^2}^\frac13||\eta||_{\htt}^\frac23\less
||M\eta||^\frac23||\eta||_{Y_0}^\frac43+||\eta||_{Y_0}^2\,.
\ee
The term $b(t)[U,\eta]$ is an indication of a  certain ``quasi-linearity" of the equation, and will be handled differently, by a suitable ``time renormalization". We  write~\rf{nstab2} as follows (keeping the same meaning of $b(t)$ as above).
\be\la{ye4}
\eta_t+(1+\ve b(t))[U,\eta+\tv]-\ve b(t)[U,\tv-b(t)U]+\ve[\tv-b(t)U,\eta]=0\,.
\ee
Note that in view of~\rf{ye2b} we have
\be\la{ye5}
||[U,\tv-b(t)U]||_{Y_0}=||[U,w]||_{Y_0}\less ||\eta||_{Y_0}\,,
\ee

Assume now again that a solution $\eta$ of~\rf{nstab2} with an initial condition $\eta_0$ for which we have~\rf{nstab2} satisfies~\rf{globalGamma}. Then
\be\la{estb1}
|b(t)|=|\tv_\theta(0,t)|=|H\eta(0,t)|\less ||\eta||_{Y_0}\less \Gamma e^{-\beta t}\,,
\ee
where we have used Lemma~\ref{hilbertonY}\,.
We rewrite the equation~\rf{ye4} as
\be\la{estb2}
\frac {\partial\eta}{(1+\ve b(t))\partial t}=\tilde L\eta+f\,,
\ee
where
\be\la{estb3}
f=\frac{\ve b(t)}{1+\ve b(t)}[U,\tv-b(t)U]-\frac \ve{1+\ve b(t)} [\tv-b(t)U,\eta]\,.
\ee
We define a ``renormalized" time variable $s$ by
\be\la{estb4}
\frac{ds}{dt}=1+\ve b(t)\,,\quad s|_{t=0}=0\,.
\ee
It is easy to see that
\be\la{estb5}
e^{-\ve\Gamma}\le \frac {e^{-\beta s}}{e^{-\beta t}}\le e^{\ve\Gamma}\,,\quad t\ge 0.
\ee
We write~\rf{estb2} as
\be\la{estb6}
\eta_s=\tL\eta+g\,,
\ee
where $g$ is defined by $g(s)=f(t)$. Collecting the various estimates, we see that
\be\la{estb7}
||g(s)||_{Y_0}\less \ve \Gamma^2 e^{-\frac43\beta t}\,.
\ee
From Duhamel's formula we see that
\be\la{estb8}
||\eta(s)||_{Y_0}=C_0e^{-\beta s} + \tilde c\ve\Gamma^2
\int_0^s C_0e^{-\beta(s-s')-\frac 43\beta s'}\,ds'=C_0(1+c\ve \Gamma^2)e^{-\beta s}\,.
\ee
Going back to the variable $t$, we have
\be\la{estb9}
||\eta(t)||_{Y_0}\le e^{c\ve\Gamma}(1+c\ve\Gamma^2)C_0e^{-\beta t}\,.
\ee
\subsection{Proof of Theorem~\ref{main}}
Let us set $\Gamma=2C_0$ and let us choose $\ve>0$ so that
\be\la{pf1}
\max\left(\,e^{c\ve\Gamma}\,\,,\,\,e^{c\ve\Gamma}(1+c\ve\Gamma^2)C_0\,\right)=e^{c\ve\Gamma}(1+c\ve\Gamma^2)C_0<\Gamma\,,
\ee
where we have used that $C_0\ge 1$. Let us consider the local solution $\eta(t)$ with $\eta(0)=\eta_0$. By continuity, the bounds~\rf{global} will be satisfied on some open time interval. If the bounds were not satisfied for all time, there would be the first moment of time $T$ when we will have equality in one of the inequalities ~\rf{global}. However, this is not possible due to the choice of $\ve$ and the bounds~\rf{estM2} and~\rf{estb9}.\cara

\subsection{Proof of Theorem~\ref{mainresult}}
Theorem~\rf{mainresult} follows from Theorem~\ref{main} by using the changes of variables detailed in the beginning of Section~\ref{nlsection}.\cara

\centerline{\bf Acknowledgement}

\smallskip
\noindent
The authors are grateful for helpful discussions with Angel Castro, Tarek Elgindi, Alex Kiselev, and Steve Preston.

The research of VS was supported in part by grants DMS 1362467 and DMS 1664297
 from the National Science Foundation. HJ was supported in part by DMS-1600779 and he is also grateful to IAS where part of the work was carried out.

\end{document}